\documentclass{article}
\usepackage{graphicx,url}
\usepackage{epsfig,epsf,listings}
\usepackage{amsmath,amssymb,amsfonts,bbm,mathrsfs}
\usepackage{caption}
\usepackage{float}
\usepackage[version=3]{mhchem}
\usepackage{ulem}
\usepackage{subfigure}
\usepackage{fancyhdr}
\usepackage{epstopdf}
\usepackage[utf8]{inputenc}
\usepackage[utf8]{inputenc}
\usepackage{enumerate}
\usepackage{mathbbol}
\usepackage{dsfont}
\usepackage{multirow,verbatim,color}
\usepackage{algorithm}
\usepackage{algpseudocode}
\usepackage{lineno,hyperref}
\usepackage{graphicx,latexsym,amsthm,times}
\usepackage{epsfig}
\usepackage{fancyhdr}
\usepackage{color,ulem}
\usepackage[english]{babel}
\usepackage{url}
\usepackage{amsfonts}
\usepackage{color}
\usepackage[margin=0.5in]{geometry}
\usepackage{amsfonts,bm}
\usepackage{upgreek}
\usepackage{setspace}
\usepackage{float}
\usepackage{tabularx}
\usepackage{xcolor}
\usepackage{graphics}
\usepackage{hyperref}
\usepackage[numbers]{natbib}
\usepackage{soul}

\newcommand{\Ns}{N_{\rm{state}}}
\newcommand{\Nt}{N_t}

\modulolinenumbers[5]

\renewcommand{\xi}{\overline{x}}

\renewcommand{\thispagestyle}[2]{}

\fancypagestyle{plain}{
        \fancyhead{}
        \fancyhead[C]{first page center header}
        \fancyfoot{}
        \fancyfoot[C]{first page center footer}
}
\begin{document}
\thispagestyle{empty}
\setcounter{page}{0}

\begin{Huge}
\begin{center}
Computer Science Technical Report CSTR-{\tt1} \\
\today
\end{center}
\end{Huge}
\vfil
\begin{huge}
\begin{center}
Azam Moosavi, Razvan Stefanescu, Adrian Sandu
\end{center}
\end{huge}

\vfil
\begin{huge}
\begin{it}
\begin{center}
``{\tt 
Multivariate predictions of local reduced-order-model errors and dimensions
}''
\end{center}
\end{it}
\end{huge}
\vfil

\begin{large}
\begin{center}
Computational Science Laboratory \\
Computer Science Department \\
Virginia Polytechnic Institute and State University \\
Blacksburg, VA 24060 \\
Phone: (540)-231-2193 \\
Fax: (540)-231-6075 \\ 
Email: \url{razvan.stefanescu@spire.com} \\
Web: \url{http://csl.cs.vt.edu}
\end{center}
\end{large}

\vspace*{1cm}

\begin{tabular}{ccc}
\includegraphics[width=2.5in]{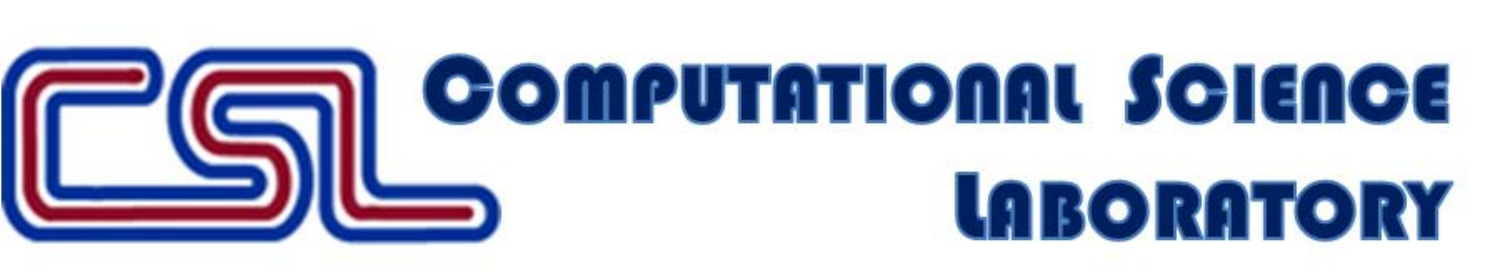}
&\hspace{2.5in}&
\includegraphics[width=2.5in]{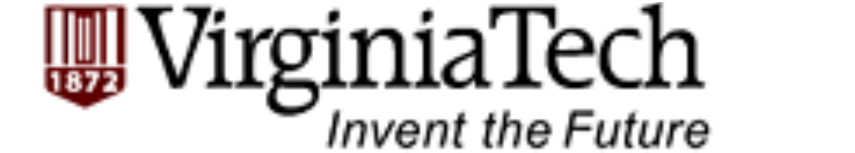} \\
{\bf\em Innovative Computational Solutions} &&\\
\end{tabular}

\newpage


\begin{abstract}
This paper introduces multivariate input-output models to predict the errors and bases dimensions of local parametric Proper Orthogonal Decomposition reduced-order models. We refer to these multivariate mappings as the MP-LROM models. We employ Gaussian Processes and Artificial Neural Networks to construct approximations of these multivariate mappings. Numerical results with a viscous Burgers model illustrate the performance and potential of the machine learning based regression MP-LROM models to approximate the characteristics of parametric local reduced-order models. The predicted reduced-order models errors are compared against the multi-fidelity correction and reduced order model error surrogates methods predictions, whereas the predicted reduced-order dimensions are tested against the standard method based on the spectrum of snapshots matrix. Since the MP-LROM models incorporate more features and elements to construct the probabilistic mappings they achieve more accurate results. However, for high-dimensional parametric spaces, the MP-LROM models might suffer from the curse of dimensionality. Scalability challenges of MP-LROM models and the feasible ways of addressing them are also discussed in this study.
\paragraph{keywords}
{Local reduced-order models, Proper Orthogonal Decomposition, regression machine learning techniques.}
\end{abstract}




\section{Introduction}
\label{sect:Intro}

Many physical phenomena are described mathematically by partial differential equations (PDEs),  and, after applying suitable discretization schemes, are simulated on a computer. PDE-based models frequently require calibration and parameter tuning in order to provide realistic simulation results. Recent developments in the field of uncertainty quantification  \cite{le2010spectral,smith2013uncertainty,grigoriu2012stochastic,cacuci2005sensitivity} provide
the necessary tools for validation of such models even in the context of variability and lack of knowledge of the input parameters. Techniques to propagate uncertainties through models include direct evaluation for linearly parametric models, sampling methods such as Monte Carlo \cite{shapiro2003monte}, Latin hypercube \cite{helton2003latin} and quasi-Monte Carlo techniques \cite{lemieux2009monte}, perturbation methods \cite{cacuci2003sensitivity,Cacuci2015687,cacuci2015second} and spectral representation \cite{le2010spectral,eldred2009comparison,alekseev2011estimation}. While stochastic Galerkin methods \cite{le2010spectral} are intrusive in nature, Monte Carlo sampling methods \cite{shapiro2003monte} and stochastic collocations \cite{eldred2009comparison} do not require the modification of existing codes and hence they are non-intrusive. While uncertainty propagation techniques can measure the impact of uncertain parameters on some quantities of interest, they often become infeasible due to the large number of model realizations requirement. Similar difficulties are encountered when solving Bayesian inference problems since sampling from posterior distribution is required.


The need for computational efficiency motivated the development of surrogate models such as response surfaces, low resolution, and reduced-order models. Data fitting or response surface models \cite{smith2013uncertainty} are data-driven models. The underlying physics remain unknown and only the input-output behavior of the model is considered. Data fitting can use techniques such as regression, interpolation, radial basis function, Gaussian Processes, Artificial Neural Networks and other supervised machine-learning methods. The latter techniques can automatically detect patterns in data, and one can use them to predict future data under uncertainty in a probabilistic framework \cite{murphy2012machine}. While easy to implement due to the non-intrusive nature, the prediction abilities may suffer since the governing physics are not specifically accounted for.


Low-fidelity models attempt to reduce the computational burden of the high-fidelity models by neglecting some of the physical aspects (e.g., replacing Navier-Stokes and Large Eddy Simulations with inviscid Euler's equations and Reynolds-Averaged Navier-Stokes \cite{gano2005hybrid,sagaut2006large,wilcox1998turbulence}, or
decreasing the spatial resolution \cite{Courtier_Thepaut1994,tremolet2007incremental}). The additional approximations, however, may considerably degrade the physical solution with only a modest decrease of the computational load.

Reduced basis \cite{porsching1985estimation,BMN2004,grepl2005posteriori,rozza2008reduced,Dihlmann_2013} and Proper Orthogonal Decomposition \cite{karhunen1946zss,loeve1955pt,hotelling1939acs,lorenz1956eof,lumley1967structure}
are two of the popular reduced-order modeling (ROM) strategies available in the literature. Data analysis is conducted to extract basis functions from experimental data or detailed simulations of high-dimensional systems (method of snapshots \cite{Sir87a, Sir87b, Sir87c}), for subsequent use in Galerkin projections that yield low dimensional dynamical models. While these type of models are physics-based and therefore require intrusive implementations, they are usually more robust than data fitting and low-fidelity models. However, since surrogate model robustness depends heavily on the problem, it must be carefully analyzed especially for large-scale nonlinear dynamical systems.

ROM robustness in a parametric setting can be achieved by constructing a global basis \cite{hinze2005proper,prud2002reliable}, but this strategy generates large dimensional bases that may lead to slow reduced-order models. Local approaches have been designed for parametric or time domains generating local bases for both the state variables \cite{Rapun_2010,dihlmann2011model} and non-linear terms \cite{eftang2012parameter,peherstorfer2014localized}. A recent survey of state-of-the-art methods in projection-based parametric model reduction is available in \cite{benner2015survey}.

In this study, we propose multivariate data fitting models to predict the local parametric Proper Orthogonal Decomposition reduced-order models errors and bases dimensions. We refer to them as MP-LROM models. Let us consider a local parametric reduced-order model of dimension $K_{POD}$ constructed using a high-fidelity solution associated with the parameter configuration $\mu_p$.

Our first MP-LROM model consists in the mapping $\{\mu, \mu_p, K_{POD}\} \mapsto \log\varepsilon_{\mu,\mu_p,K_{POD}}^{HF}$, where $\varepsilon_{\mu,\mu_p,K_{POD}}^{HF}$ is the error of the local reduced-order model solution with respect to the high-fidelity solution for a viscosity parameter configuration $\mu$. Our proposed approach is inspired from the multi-fidelity correction (MFC) \cite{alexandrov2001approximation} and reduced order model error surrogates method (ROMES) \cite{drohmann2015romes}.  MFC \cite{alexandrov2001approximation,eldred2004second,gano2005hybrid,huang2006sequential}
has been developed for low-fidelity models in the context of optimization. The MFC model simulates the input-output relation $\mu \mapsto \varepsilon_{\mu}^{HF}$, where $\varepsilon_{\mu}^{HF}$ is the low-fidelity model error  depending on a global reduced basis with a constant reduced-order model dimension. The ROMES method \cite{drohmann2015romes} introduced the concept of error indicators for global reduced-order models and generalized the MFC framework by approximating the mapping $\rho(\mu) \mapsto \log\varepsilon_{\mu}^{HF}$. The error indicators $\rho(\mu)$ include rigorous error bounds and reduced-order residual norms. No variation of the reduced basis dimension was taken into account. By estimating the log of the reduced-order model error instead of the error itself, the input-output map exhibits a lower variance as shown by our numerical experiments as well as those in \cite{drohmann2015romes}.

The second proposed MP-LROM model addresses the issue of a-priori selection of the reduced basis dimension for a prescribed accuracy of the reduced solution. The standard approach is to analyze the spectrum of the snapshots matrix, and use the largest singular value removed from the expansion to estimate the accuracy level \cite{volkwein2007proper}. To also take into account the error due to the full-order-model equations projection in the reduced space, here we propose the mapping  $\{\mu_p, \log\varepsilon_{\mu_p,\mu_p,K_{POD}}^{HF}\} \mapsto K_{POD}$ to predict the dimension of a local parametric reduced-order model given a prescribed error threshold.

To approximate the mappings $\{\mu, \mu_p, K_{POD}\} \mapsto \log\varepsilon_{\mu,\mu_p,K_{POD}}^{HF}$ and $\{\mu_p, \log\varepsilon_{\mu_p,\mu_p,K_{POD}}^{HF}\} \mapsto K_{POD}$, we propose regression models constructed using Gaussian Processes (GP) \cite{slonski2011bayesian,lilley2004gaussian} and Artificial Neural Networks (ANN). In the case of one dimensional Burgers model, the resulted MP-LROM error models are accurate and their predictions are compared against those obtained by the MFC and ROMES models. The predicted dimensions of local reduced-order models using our proposed MP-LROM models are more accurate than those derived using the standard method based on the spectrum of snapshots matrix.

The remainder of the paper is organized as follows. Section \ref{sect:ROM} reviews the reduced-order modeling parametric framework. The MP-LROM models and the regression machine learning methods used in this study to approximate the MP-LROM mappings are described in details in Section \ref{sect:MP-LROM}.
Section \ref{sect:experm} describes the viscous 1D-Burgers model and compares the performances of the MP-LROM and state of the art models.
Conclusions are drawn in Section \ref{sect:conc}.

\section{Parametric reduced-order modeling}
\label{sec:ROM}
\label{sect:ROM}

Proper Orthogonal Decomposition has been successfully applied in numerous applications such as compressible flow \citep{Rowley2004} and computational fluid dynamics \citep{Kunisch_Volkwein_POD2002,Rowley2005,Willcox02balancedmodel}, to mention a few. It can be thought of as a Galerkin approximation in the state variable built from functions corresponding to the solution of the physical system at specified time instances. A system reduction strategy for Galerkin models of fluid flows based on a partition in slow, dominant, and fast modes, has been proposed in \cite{Noack2010}. Closure models and stabilization strategies for POD of turbulent flows have been investigated in \cite{San_Iliescu2013,wells2015regularized}.

In this paper we consider discrete inner products (Euclidean dot product), though continuous products may be employed as well.
Generally, an unsteady problem can be
written in semi-discrete form as an initial value problem; i.e., as a system of nonlinear ordinary differential equations
\begin{equation}
\label{eqn::-4}
\frac{d{\bf x}(\mu,t)}{dt} = {\bf F}({\bf x},t,\mu),~~~~{\bf x}(\mu,0) = {\bf x}_0 \in \mathbb{R}^{\Ns},
\quad \mu \in \mathcal{P}.
\end{equation}
The input-parameter $\mu$ typically characterizes the physical properties of the flow. By $\mathcal{P}$ we denote the input-parameter space. For a given parameter configuration $\mu_p$ we select an ensemble of $N_t$ time instances of the flow
${\bf x}(\mu_p, {t_1}),\ldots,{\bf x}(\mu_p , t_{N_t}) \in \mathbb{R}^{\Ns}$, where ${\Ns}$ is the total number of discrete model variables, and $N_t \in \mathbb{N^*}$.  The POD method chooses an orthonormal basis
$U_{\mu_p}=[{\bf u}_{1}^{\mu_p} ~~\cdots~~ {\bf u}_{K_{POD}}^{\mu_p}] \in \mathbb{R}^{{\Ns}\times K_{POD}}$, such that the mean square error between ${\bf x}(\mu_p,t_i)$ and the POD expansion
 ${\bf x}^\textsc{pod}_{\mu_p}(t_i) = U_{\mu_p}{\bf \tilde x_{\mu_p}}(\mu,t_i)$, ${\bf \tilde  x_{\mu_p}}(\mu,t_i)= U_{\mu_p}^T {\bf x}(\mu_p , t_i)
  \in \mathbb{R}^ {K_{POD}} $, is minimized on average. The POD space dimension $K_{POD} \ll {\Ns}$ is appropriately chosen to capture the dynamics of the flow. Algorithm \ref{euclid} describes the reduced-order basis construction procedure \cite{stefanescu2014comparison}.

\begin{algorithm}
 \begin{algorithmic}[1]
 \State Compute the singular value decomposition for the snapshots matrix $ [{\bf x}(\mu_p, {t_1})~ \cdots ~{\bf x}(\mu_p, {t_{N_t}})]= \bar U_{\mu_p} \Sigma_{\mu_p} {\bar V}^T_{\mu_p},$ with the singular vectors matrix $\bar U_{\mu_p} =[{\bf u}_1^{\mu_p}~~ \cdots ~~{\bf u}_{N_t}^{\mu_p}].$
 \State Using the singular-values $\lambda_1\geq \lambda_2\geq  \ldots \geq \lambda_{N_t} \geq 0$ stored in the diagonal matrix $\Sigma_{\mu_p}$, define $I(m)= {\sum_{i=1}^m \lambda_i^2/(\sum_{i=1}^{t_{N_t}} \lambda_i^2})$.
\State Choose $K_{POD}$, the dimension of the POD basis, such that $ K_{POD}={\rm arg}\min_m \{I(m):I(m)\geq \gamma\}$ where $0 \leq \gamma \leq 1$ is the percentage of total information captured by the reduced space $\mathcal{X}^{K_{POD}}=\textnormal{range}(U_{\mu_p})$. It is common to select $\gamma=0.99$. The basis $U_{\mu_p}$ consists of the first $K_{POD}$ columns of $\bar U_{\mu_p}$.
 \end{algorithmic}
 \caption{POD basis construction}
 \label{euclid}
\end{algorithm}

Next, a Galerkin projection of the full model state \eqref{eqn::-4} onto the space $\mathcal{X}^{K_{POD}}$ spanned by the POD basis elements is used to obtain the reduced-order model
\begin{equation}\label{eqn::-3}
 \frac{d{\bf \tilde x}_{\mu_p}(\mu,t)}{dt} = U_{ \mu_p}^T\,{\bf F}\bigg(U_{\mu_p}{\bf \tilde x}_{\mu_p}(\mu,t), t, \mu \bigg),
 \quad {\bf \tilde x}_{\mu_p}(\mu,0)= U_{\mu_p}^T\,{\bf x}_0.
 \end{equation}
 The notation ${\bf \tilde x}_{\mu_p}(\mu,t)$ expresses the solution dependence on the varying parameter $\mu$ and also on $\mu_p$ the configuration whose associated high-fidelity trajectory was employed to generate the POD basis. While being accurate for $\mu=\mu_p$, the reduced model  \eqref{eqn::-3} may lose accuracy when moving away from the initial setting. Several strategies have been proposed to derive a basis that spans the entire parameter space. These include the reduced basis method combined with the use of error estimates \cite{rozza2008reduced,quarteroni2011certified,prud2002reliable}, global POD \cite{taylor2004towards,schmit2003improvements}, Krylov-based sampling methods \cite{daniel2004multiparameter,weile1999method}, and  greedy techniques  \cite{haasdonk2008reduced,nguyen2009reduced}. The fundamental assumption used by these approaches is that a smooth low-dimensional global manifold characterizes the model solutions over the entire parameter domain. The purpose of our paper is to estimate the solution error and dimension of the reduced-order model \eqref{eqn::-3} that can be subsequently used to generate a global basis for the parameter space.

\section{Multivariate prediction of local reduced-order models characteristics (MP-LROM) \label{sect:MP-LROM}}

We propose multivariate input-output models
 \begin{equation}\label{eqn:general_MP-LROM}
 \phi: {\bf z} \mapsto {{y}},
\end{equation}
${\bf z} \in \mathbb{R}^r$, to predict characteristics $y \in \mathbb{R}$ of local parametric reduced-order models \eqref{eqn::-3}.


\subsection{Error Model} \label{sect:error_MP-LROM}

Inspired from the MFC and ROMES methodologies we introduce an input-output model to predict the level of error $\varepsilon_{\mu,\mu_{p},K_{POD}}^{HF}$, where
\begin{equation}\label{eqn:level_error_ML_ROM}
\begin{array}{lr}
 {\varepsilon_{\mu,\mu_{p},K_{POD}}^{HF}} = \\
\| {\bf x}(\mu,t_1) - U_{ \mu_{p}}{\bf \tilde x}_{\mu_{p}}(\mu,t_1) \quad {\bf x}(\mu,t_2) - U_{ \mu_{p}}{\bf \tilde x}_{\mu_{p}}(\mu,t_2) \quad \cdots \quad {\bf x}(\mu,t_{N_t}) - U_{ \mu_{p}}{\bf \tilde x}_{\mu_{p}}(\mu,t_{N_t})   \|_F.
\end{array}
\end{equation}
Here $\|\cdot\|_F$ denotes the Frobenius norm, and $K_{POD}$ is the dimension of the reduced-order model. In contrast with ROMES and MFC models that predict the error of global reduced-order models with fixed dimensions, using univariate functions, here we propose a multivariate model
 \begin{equation}\label{eqn:MP-LROM-error}
 \phi_{MP-LROM}^e: \{\mu, \mu_{p}, K_{POD}\} \mapsto \log \varepsilon_{\mu,\mu_{p},K_{POD}}^{HF}
\end{equation}
to predict the error of local parametric reduced-order models \eqref{eqn::-3} of various dimensions. Since the dimension of basis usually influences the level of error we include it among the input variables. To design models with reduced variances we look to approximate the logarithm of the error as suggested in \cite{drohmann2015romes}.

For high-dimensional parametric spaces, ROMES method handles well the curse of dimensionality with their proposing univariate models. In combination with active subspace method \cite{constantine2014active}, we can reduce the number of input variables in case the amount of variability in the parametric space is mild. This will increase our error model feasibility even for high-dimensional parametric space.

\subsection{Dimension of the reduced basis}

The basis dimension represents one of the most important characteristic of a reduced-order model. The reduced manifold dimension directly affects both the on-line computational complexity of the reduced-order model and its accuracy \cite{kunisch2001galerkin,Hinze_Wolkwein2008,fahl2003reduced}. {By increasing the dimension of the basis, the projection error usually decreases and the accuracy of the reduced-order model is enhanced. However this is not necessarily valid as seen in \cite[Section 5]{rowley2004model}.} Nevertheless the spectrum of the snapshots matrix offers guidance regarding the choice of the reduced basis dimension when some prescribed reduced-order model error is desired. However the accuracy depends also on the `in-plane' error, which is due to the fact that the full-order-model equations are projected on the reduced subspace \cite{MRathinam_LPetzold_2003a,homescu2005error}.

We seek to predict the dimension of the local parametric reduced-order model \eqref{eqn::-3} by accounting for both the orthogonal projection error onto the subspace, which is computable by the sum of squares of singular values, and the `in-plane' error. As such we propose to model the mapping
 \begin{equation}\label{eqn:eqn:MP-LROM-dimension}
 \phi_{MP-LROM}^d: \{\mu_{p}, \log \varepsilon_{\mu_{p},\mu_{p},K_{POD}}^{HF}\} \mapsto K_{POD}.
\end{equation}
Once such model is available, given a positive threshold  $\bar{\varepsilon}$ and a parametric configuration $\mu_p$, we will be able to predict the dimension $K_{POD}$ of the basis $U_{\mu_p}$, such that the reduced-order model error satisfies
\begin{equation}
\label{eqn:level_error3}
\| {\bf x}(\mu_p,t_1) - U_{ \mu_p}{\bf \tilde x}_{\mu_p}(\mu_p,t_1) \quad {\bf x}(\mu_p,t_2) - U_{ \mu_p}{\bf \tilde x}_{\mu_p}(\mu_p,t_2) \quad \cdots \quad {\bf x}(\mu_p,t_{N_t}) - U_{ \mu_p}{\bf \tilde x}_{\mu_p}(\mu_p,t_{N_t})   \|_F
\approx \bar{\varepsilon}.
\end{equation}
%

\subsection{Supervised Machine Learning Techniques}
\label{sect:prob_fram}

In order to estimate the level of reduced-order model solution error ${\varepsilon_{\mu,\mu_{p_j},K_{POD}}^{HF}}$ \eqref{eqn:level_error_ML_ROM} and the reduced basis dimension $K_{POD}$, we will use regression machine learning methods to approximate the maps $\phi_{MP-LROM}^e$ and $\phi_{MP-LROM}^d$ described in \eqref{eqn:MP-LROM-error} and \eqref{eqn:eqn:MP-LROM-dimension}.

Artificial Neural Networks and Gaussian Processes are used to build a probabilistic model $ \phi: {\bf z} \mapsto \hat{y} $, where $\phi$ is a transformation function that learns through the input features ${\bf z}$  to estimate the deterministic output $y$ \cite{murphy2012machine}. As such, these probabilistic models are approximations of the mappings introduced in \eqref{eqn:general_MP-LROM}. The input features ${\bf z}$ can be either  categorical or ordinal. The real-valued random variable $\hat{y}$ is expected to have a low variance and reduced bias. The features of $ {\bf z} $ should be descriptive of the underlying problem at hand \cite{bishop2006pattern}. The accuracy and stability of estimations are assessed using the K-fold cross-validation technique. The samples are split into K subsets (``folds''), where typically $3 \le K \le 10$. The model is trained on $K-1$ sets and tested on the $K$-th set in a round-robin fashion \cite{murphy2012machine}.  Each fold induces a specific error quantified as the average of the absolute values of the differences between the predicted and the $K$-th set values
\begin{subequations}
\begin{equation}
\label{eqn:err_fold}
\textnormal{E}_{\rm fold}=\frac{\sum_{i=1}^N | \hat{y}^i-y^i | }{N} , \quad
\textnormal{VAR}_{\rm fold}=\frac{\sum_{i=1}^N \left( \hat{y}^i - \textnormal{E}_{\rm fold} \right)^2}{N-1}, \quad \rm fold=1,2, \ldots, K,
\end{equation}
where $N$ is the number of test samples in the fold.  The error is then averaged over all folds:
\begin{equation}
\label{eqn:err_fold_average}
\textnormal{E}=\frac{\sum_{\textnormal{fold}=1}^K\,  \textnormal{E}_{\rm fold}  }{K}, \quad
\textnormal{VAR}=\frac{\sum_{\textnormal{fold}=1}^K \left(\textnormal{E}_{\rm fold}- \textnormal{E} \right)^2}{K-1}.
\end{equation}
\end{subequations}

The variance of the prediction results \eqref{eqn:err_fold} accounts for the sensitivity of the model to the particular choice of data set. It quantifies the stability of the model in response to the new training samples. A smaller variance indicates more stable predictions, however, this sometimes translates into a larger bias of the model. Models with small variance and high bias make strong assumptions about the data and tend to underfit the truth, while models with high variance and low bias tend to overfit the truth \cite{biasVar_NG} . The trade-off between bias and variance in learning algorithms is usually controlled via techniques such as regularization or bagging and boosting \cite{bishop2006pattern}.

In what follows we briefly review the Gaussian Process and Artificial Neural Networks techniques.


\subsubsection{Gaussian process kernel method}
\label{sect:GP}
A Gaussian process is a collection of random variables, any finite number of which have a joint Gaussian distribution  \cite{rasmussen2006gaussian}. A Gaussian process is fully described by its mean and covariance functions

\begin{equation}
\label{GP_Dist}
\phi(\mathbf{z}) \sim \textnormal{gp}\, \bigl({\it m}(\mathbf{z}), \mathcal{\bf K}) \bigr),
\end{equation}
where $ {\it m}(\mathbf{z})=\mathbb{E}\left[ \phi(\mathbf{z}) \right], $
and ${\bf K}$  is the covariance matrix with entries $ {K}_{i,j} = \mathbb{E} \left[\left(\phi(\mathbf{z}^i)-{\it m}(\mathbf{z}^i)\right) \left( \phi(\mathbf{z}^j)- {\it m} (\mathbf{z}^j) \right)  \right]$
\cite{rasmussen2006gaussian}.

In this work we employ the commonly used squared-exponential-covariance Gaussian kernel with

\begin{equation}
\label{eq_cov}
k:\mathbb{R}^r \times \mathbb{R}^r \rightarrow \mathbb{R},~k(\mathbf{z}^i,\mathbf{z}^j) =\sigma^2_\phi\, \exp \left(-\frac{ \left\lVert \mathbf{z}^i -  \mathbf{z}^j  \right\rVert}{2\, \hslash ^2} \  \right)+ \sigma^2_n  \, \delta_{i,j},
\end{equation}
and ${K}_{ij} = k(\mathbf{z}^i,\mathbf{z}^j)$ \cite{rasmussen2006gaussian}, where $\mathbf{z}^i $ and $\mathbf{z}^j$ are the pairs of data points in training or test samples, $\delta $ is the Kronecker delta symbol and $\|\cdot \|$ is some appropriate norm. The model \eqref{eq_cov} has three hyper-parameters. The length-scale $\hslash$ governs the correlation among data points. The signal variance $\sigma^2 _\phi \in \mathbb{R}$ and the noise variance $\sigma^2 _n \in \mathbb{R} $ govern the precision of variance and noise, respectively.

Consider a set of training data points ${\bf Z} = [\mathbf{z}^1~\mathbf{z}^2 ~~\cdots~~  \mathbf{z}^n] \in \mathbb{R}^{r \times n} $ and the corresponding noisy observations ${\bf y} = [y^1~y^2 ~~\cdots~~ y^n] \in \mathbb{R}^{1 \times n},$
\begin{equation}
\label{GP_training}
y^i=\phi(\mathbf{z}^i)+ \epsilon_i ,\quad \epsilon_i \sim \mathcal{N} \left(0, \sigma^2_n \right), \quad i = 1,\dots,n.
\end{equation}
%

Consider also the set of test points ${\bf Z}^* = [\mathbf{z}^{*1}~\mathbf{z}^{*2} ~~\cdots~~  \mathbf{z}^{*m}] \in \mathbb{R}^{r \times m}$ and the predictions ${\bf \hat{y}} = [\hat{y}^1~\hat{y}^2~~\cdots~~ \hat{y}^m] \in \mathbb{R}^{1 \times m}$,
\begin{equation}
\label{GP_test}
\hat{y}^i=\phi\left(\mathbf{z}^{*i}\right), \quad i = 1,\dots,m.
\end{equation}

For a Gaussian prior the joint distribution of training outputs ${\bf y}$ and test outputs ${\bf \hat{y}}$ is
\begin{equation}
\label{GP_prior}
\begin{bmatrix}
{\bf y}^T\\
{\bf \hat{y}^T}
\end{bmatrix}
\sim
\mathcal{N}
\left(
\begin{bmatrix} { {\bf m}}(\mathbf{Z})^T \\ { {\bf m}}(\mathbf{Z}^*)^T \end{bmatrix}\, , \,
\begin{bmatrix}
{\bf K}& {\bf K}^*\\
{\bf K}^{*T} & {\bf K}^{**}\\
\end{bmatrix}
\right),
\end{equation}
where

$${\bf m}({\bf Z}) = [{\it m}({\bf z}^1)~{\it m}({\bf z}^2) ~~\cdots~~ {\it m}({\bf z}^n)]\in \mathbb{R}^{1 \times n},~{\bf m}({\bf Z}^*)= [{\it m}({\bf z}^{*1}) ~{\it m}({\bf z}^{*2}) ~~\cdots~~ {\it m}({\bf z}^{*m})] \in \mathbb{R}^{1 \times m},$$
$${\bf K}^* = ({K_{ij}^*})_{i=1,\ldots,n;~j=1,\ldots,m} = k({\bf z}^i,{\bf z}^{j*}) \textrm{ and } {\bf K}^{**} = ({K_{ij}^{**}})_{i=1,\ldots,m;~j=1,\ldots,m} = k({\bf z}^{i*},{\bf z}^{j*}).$$

The predictive distribution represents the posterior after observing the data \cite{bishop2006pattern} and is given by

\begin{equation}
\label{GP_posterior}
p\left({\bf \hat{y}}|\mathbf{Z},{\bf y},\mathbf{Z}^* \right) \sim \mathcal{N}
\left(\, {\bf K}^{*T}{\bf K}^{-1}{\bf y}\, , \,
{\bf K}^{**}- {\bf K}^{*T} {\bf K}^{-1} {\bf K}^*\, \right),
\end{equation}
where superscript $T$ denotes the transpose operation.

The prediction of Gaussian process will depend on the choice of the mean and covariance functions, and on their hyper parameters $\hslash$, $\sigma^2 _\phi $ and $ \sigma^2_n $  which  can be inferred from the data
 $${\bm \theta}^* = [\hslash, \sigma^2 _\phi, \sigma^2_n ] = \arg\min_{{\bm \theta}}\, L({\bm  \theta}),$$
by minimizing the marginal negative log-likelihood function
\[
L({\bm \theta}) = - \log\, p({\bf y}|\mathbf{Z},{\bm \theta})=\frac{1}{2} \log \det({\mathbf{K}})
+ \frac{1}{2}  ({\bf y}-{\bf m}(\mathbf{Z}))\, {\mathbf{K}}^{-1}\, ({\bf y}-{\bf m}(\mathbf{Z}))^T + \frac{n}{2}\, \log \left( 2 \pi \right).
\]
%

\subsubsection{Artificial Neural Networks}
\label{sect:NN}
The study of Artificial Neural Networks begins in the 1910s in order to imitate human brain's biological structure.
Pioneering work was carried out by Rosenblatt, who proposed a three-layered network structure, the perceptron \cite{hagan2014neural} .
%
ANN detect the pattern of data by discovering the input--output relationships. Applications include the approximation of functions, regression analysis, time series prediction, pattern recognition, and speech synthesis and recognition \cite{jang1997neuro,ayanzadeh2011fossil}.

ANN consist of neurons and connections between the neurons (weights).
Neurons are organized in layers, where at least three layers of neurons (an input layer, a hidden layer, and an output layer) are required for construction of a neural network.

The input layer distributes input signals $\mathbf{z} = [{z}_1~{z}_2 ~~\cdots~~ {z}_r]$ to the first hidden layer. For a neural network with $L$ hidden layers and $m^{\ell}$ neurons in each hidden layer, let $ {\bf \hat{y}}^{\ell}= [\hat{y}^{\ell}_1~\hat{y}^{\ell}_2~~\cdots~~ \hat{y}^{\ell}_{m^{\ell}} ]$  be the vector of outputs from layer $\ell$,
$\mathbf{b}^\ell = [b^{\ell}_1~b^{\ell}_2~~\cdots  ~~ b^{\ell}_{m^{\ell}}]$ the biases at layer $\ell$, and ${\bf w}_{j}^\ell = [{w}_{j_1}^\ell {w}_{j_2}^\ell ~~\cdots~~ w^l_{j_{m^l}}]$ the weights connecting the neuron $j$ to the input of that layer (output of previous layer).
The vectors ${\bf \hat{y}}^{\ell}$ and ${\bf w}_{j}^\ell$ share the same dimension which varies along the layers depending on the number of input features, neurons and outputs. Then the feed-forward operation is
\[
 \begin{array}{lr}
   {x}_j^{\ell+1}={{\bf w}_j^{\ell +1}}^T {\bf \hat{y}^{\ell}} + b_j^{\ell+1} ,  \quad {\bf \hat{y}}^0= \mathbf{z},\quad  j=1,\ldots {  , }m^{\ell}{.}\\
   \hat{y}_j^{\ell+1}=\varphi \left(\mathbf{x}^{\ell+1}  \right), \quad  \ell=0, 1, \ldots, L-1.
 \end{array}
\]
{All products of previous layer output with current layer neuron weights will be summed and the bias value of each neuron will be added to obtain the vector $\mathbf{x}^\ell = [x^{\ell}_1~x^{\ell}_2~~\cdots  ~~ x^{\ell}_{m^{\ell}}]$ .} Then the final output of each layer will be obtained by passing the vector $\mathbf{x}^\ell$ through the transfer function $\varphi$, which is a
differentiable function and can be log-sigmoid, hyperbolic tangent sigmoid, or linear transfer function.
%

 The training process of ANN adjusts the weights and the biases in order to reproduce the desired outputs when fed the given inputs. The training process via the back propagation algorithm \cite{rumelhart1985learning} uses a gradient descent method to modify weights and thresholds such that the error between the desired output and the output signal of the network is minimized \cite{funahashi1989approximate}.  In supervised learning the network is provided with samples from which it discovers the relations of inputs and outputs. The output of the network is compared with the desired output, and  the error is back-propagated through the network and the weights will be adjusted. This process is repeated during several
iterations, until the network output is close to the desired output \cite{haykin2009neural}.

\section{Numerical experiments}
\label{sect:experm}

We illustrate the application of the proposed MP-LROM models to predict the error and dimension of the local parametric reduced-order models for a one-dimensional Burgers model. The 1D-Burgers model proposed herein is characterized by the viscosity coefficient. To assess the performance of the MP-LROM models constructed using Gaussian Process and Artificial Neural Networks, we employ various cross-validation tests. The dimensions of the training and testing data sets are chosen empirically based on the number of samples. For Artificial Neural Networks models the number of hidden layers and neurons in each hidden layer vary for each type of problems under study. The squared-exponential-covariance kernel \eqref{eq_cov} is used for Gaussian Process models.

The approximated MP-LROM error models  are compared against the ROMES and multi-fidelity correction models, whereas the MP-LROM models that predict the dimension of the reduced-order models are verified against the standard approach based on the spectrum of snapshots matrix.


\subsection{One-dimensional Burgers' equation}\label{subsec:Burgers}
\ifx
Here we propose two alternative approaches to select the reduced basis size that accounts for specified accuracy levels in the reduced-order model solutions.
Assume we have a probability space $(\Omega, \mathcal{F},\mathcal{P})$. These techniques employ construction of probabilistic models via ANN and GP, $ \phi: X \rightarrow \hat{y} $ where $\phi$ is the transformation function that learns through the input features $X$ to estimate the deterministic
output $y \in \mathbb{R}$ through a real-valued random variable $\hat{y} : \Omega \mapsto \mathbb{R}$.
\fi

Burgers' equation is an important partial differential equation from fluid mechanics \cite{burgers1948mathematical}. The evolution of the velocity $u$ of a fluid evolves according to
\begin{equation}
\frac{\partial u}{\partial t} + u\frac{\partial u}{\partial x} = \mu \frac{\partial^2 u}{\partial x^2}, \quad x \in [0,L], \quad t \in (0,t_\textnormal{f}],\label{eqn:Burgers-pde}
\end{equation}
with $t_\textnormal{f} = 1$ and $L=1$. Here $\mu$ is the viscosity coefficient. 

\begin{figure}[h]
  \centering
\includegraphics[scale=0.37]{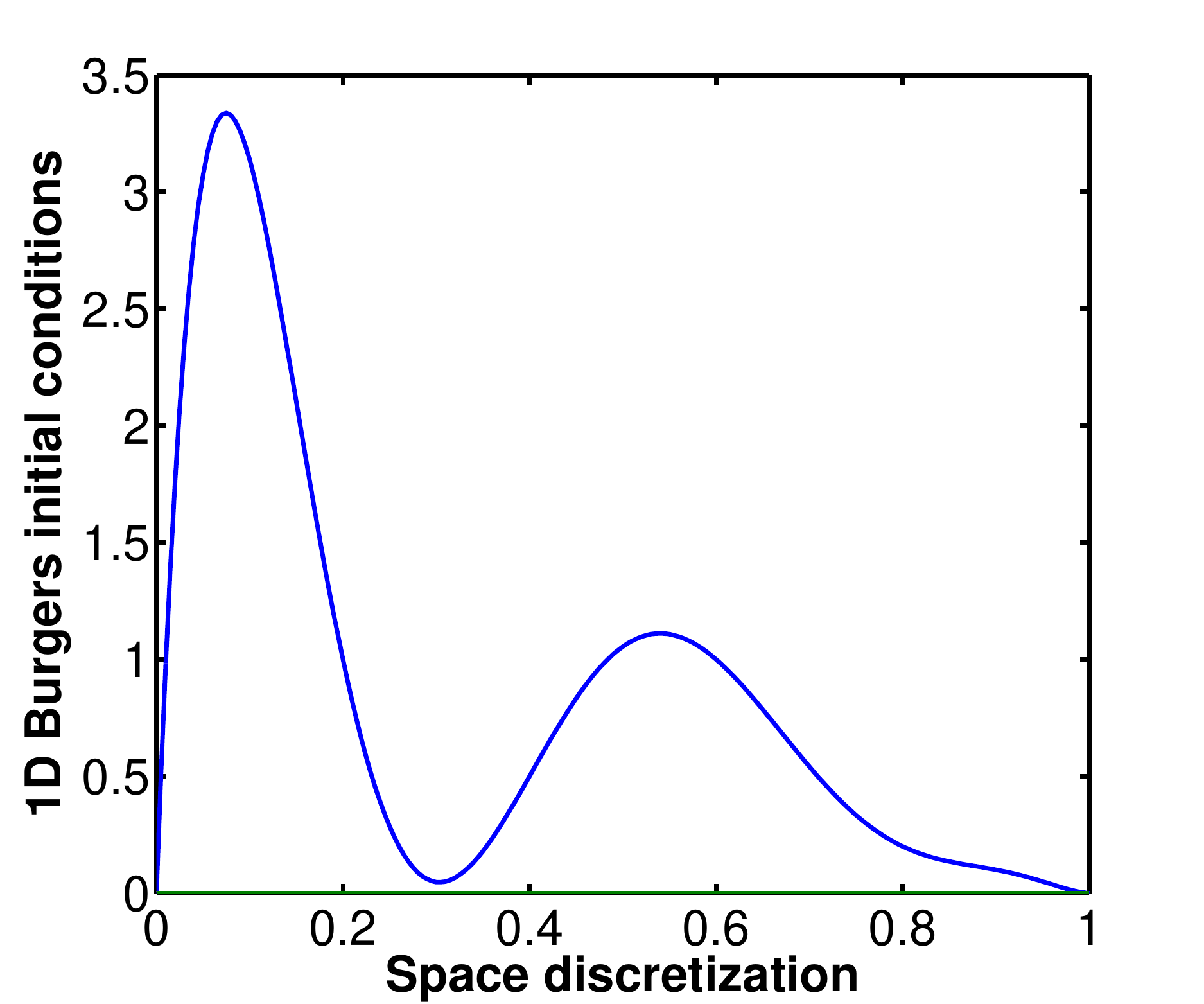}
\caption{Seventh order polynomial used as initial conditions for 1D Burgers model.
\label{Fig::1D-Burgers-IC}}
\end{figure}

 The model has homogeneous Dirichlet boundary conditions $u(0,t) = u(L,t) = 0$, $t \in (0,t_\textnormal{f}]$. For the initial conditions, we used a seventh order polynomial constructed using the least-square method and the data set $\{(0,0);~(0.2,1);~(0.4,0.5);~(0.6,1);
 ~(0.8,0.2);\newline~~(0.9,0.1);~(0.95,0.05);~(1,0) \}$. We employed the polyfit function in Matlab and the polynomial is shown in Figure \ref{Fig::1D-Burgers-IC}.


The discretization uses a spatial mesh of $N_s$ equidistant points on $[0,L]$, with $\Delta x=L/(N_s-1)$. A uniform temporal mesh with $N_t$ points  covers the interval $[0,t_\textnormal{f}]$, with $\Delta t=t_\textnormal{f}/(N_t-1)$. The discrete velocity vector is ${\boldsymbol u}(t_j)\approx [u(x_i,t_j)]_{i=1,2, \ldots,\Ns} \in \mathbb{R}^{\Ns}$,  $j=1,2, \ldots, N_t$, where $\Ns=N_s-2$  (the known boundaries are removed). The semi-discrete version of the model \eqref{eqn:Burgers-pde} is

\begin{equation}\label{eqn:Burgers-sd}
 {\bf u}'  =  -{\bf u}\odot A_x{\boldsymbol u} + \mu A_{xx}{\boldsymbol u},
\end{equation}
where ${\bf u}'$ is the time derivative of ${\bf u}$, and $A_x,A_{xx}\in \mathbb{R}^{\Ns\times \Ns}$ are the central difference first-order and second-order space derivative operators, respectively, which take into account the boundary conditions, too.  The model is implemented in Matlab and the backward Euler method is employed for time discretization. The nonlinear algebraic systems are solved using the Newton-Raphson method and the allowed number of Newton iterations per each time step is set to $50$. The solution is considered to have converged when the Euclidean norm of the residual is less then $10^{-10}$.

The viscosity parameter space $\mathcal{P}$ is set to the interval $[0.01,1]$. Smaller values of $\mu$ correspond to sharper gradients in the solution, and lead to dynamics more difficult to accurately approximate using reduced-order models.

\begin{figure}[h]
  \centering
\includegraphics[scale=0.37]{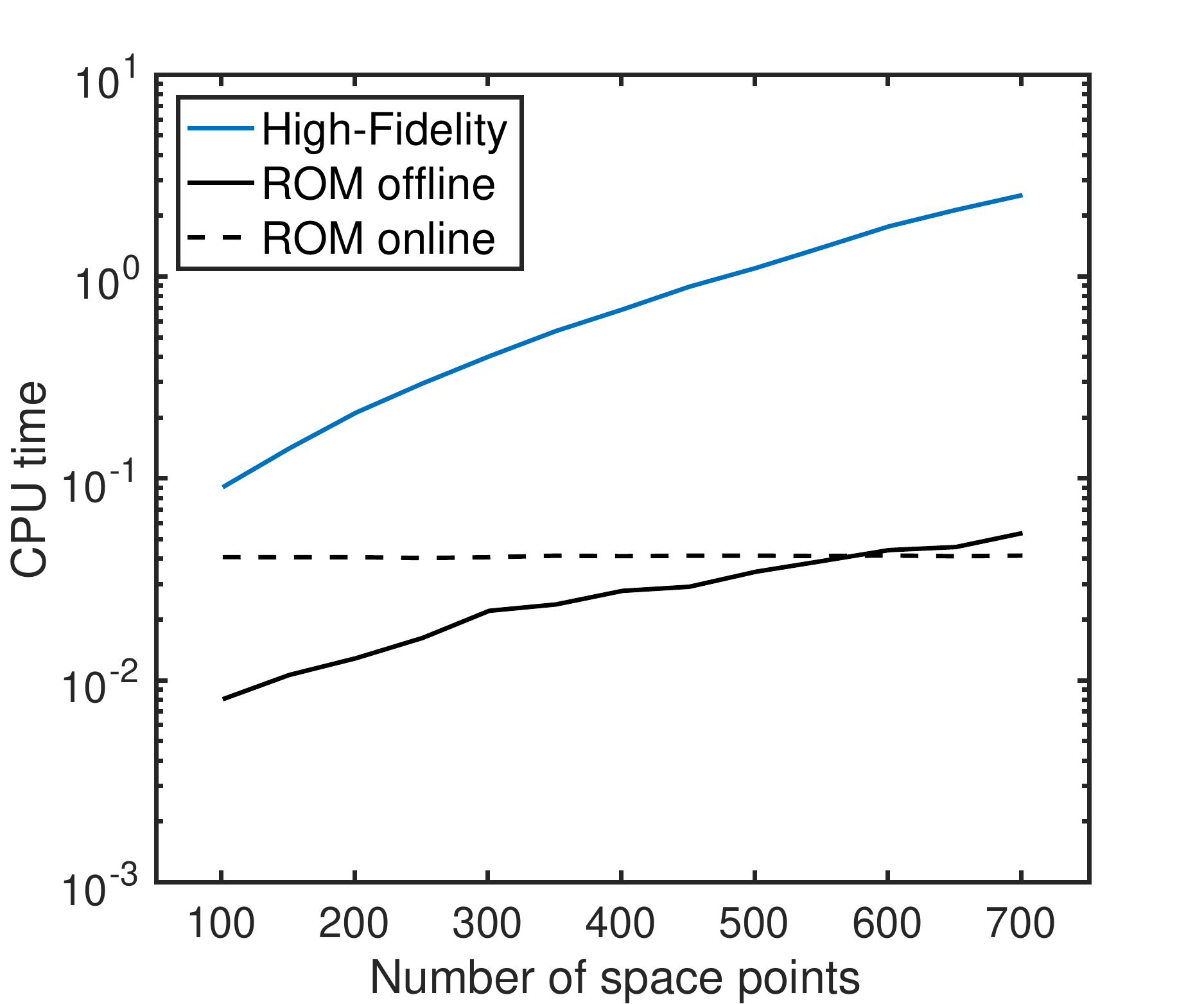}
\caption{Computational efficiency of the tensorial POD 1D Burgers model. CPU time is given in seconds.
\label{Fig::1D-Burgers-CPU_time}}
\end{figure}

The reduced-order models are constructed using POD method whereas the quadratic nonlinearities are computed via tensorial POD \cite{stefanescu2014comparison} for efficiency. A floating point operations analysis of tensorial POD, POD and POD/DEIM for $p^{\textrm{th}}$ order polynomial nonlinearities is available in \cite{stefanescu2014comparison}.
The computational efficiency of the tensorial POD 1D Burgers model can be noticed in Figure \ref{Fig::1D-Burgers-CPU_time}. Both on-line and off-line computational costs are shown. Here we selected $\mu = \mu_p = 0.7$, $\Nt = 301,$ POD dimension $K_{POD} = 9$, and we let the number of space points $N_s$ to vary. For $N_s = 201$ and $701$, the tensorial POD model is $5.17 \times$ and $61.12 \times$ times faster than the high-fidelity version. The rest of our numerical experiments uses $N_s = 201$ and $\Nt = 301$.

\subsection{Multivariate prediction of local reduced-order models characteristics (MP-LROM) using regression machine learning methods}

\subsubsection{Error estimation of local ROM solutions }
\label{sect:err_estimate}

Here, we will use GP and ANN to approximate the MP-LROM error model introduced in \eqref{eqn:MP-LROM-error}. The approximated models have the following form

 \begin{equation}\label{eqn:prob_model_scale}
  \phi_{MP-LROM}^e: \{\mu,\mu_p,K_{POD}\} \mapsto \widehat{\log{{\varepsilon}}}_{\mu,\mu_p,K_{POD}}^{HF},
\end{equation}
where the input features include a viscosity parameter value $\mu$, a parameter value $\mu_p$ associated with the full model run that generated the basis $U_{\mu_p}$, and the dimension of the reduced manifold $K_{POD}$. The target is the estimated logarithm of error of the reduced-order model solution at $\mu$ using the basis $U_{\mu_p}$ and the corresponding reduced operators computed using the Frobenius norm

\begin{equation} \label{eqn:param_rang_err}
\begin{array}{lr}
\log {{\varepsilon}}_{\mu,\mu_p,K_{POD}}^{HF} =\\
 \log\Bigg(\| {\bf x}(\mu,t_1) - U_{ \mu_p}{\bf \tilde x}_{\mu_p}(\mu,t_1) \quad {\bf x}(\mu,t_2) - U_{ \mu_p}{\bf \tilde x}_{\mu_p}(\mu,t_2) \quad \cdots \quad {\bf x}(\mu,t_{N_t}) - U_{ \mu_p}{\bf \tilde x}_{\mu_p}(\mu,t_{N_t})   \|_F\Bigg).
\end{array}
\end{equation}

The probabilistic models described generically in equation \eqref{eqn:prob_model_scale} are just approximations of the MP-LROM model \eqref{eqn:MP-LROM-error} and have errors. For our experiments, the data set includes $10$ and $100$ equally distributed values of $\mu_p$ and $\mu$ over the entire parameter region; i.e., $\, \mu_p \in \{0.1, 0.2,\ldots, 1\}$ and $\, \mu \in \{ 0.01, \ldots, 1\}$, $12$ reduced basis dimensions $K_{POD}$ spanning the interval $\{4,5,\ldots,14,15\}$ and the reduced-order model logarithm of errors $\log{{\varepsilon}}_{\mu,\mu_p,K_{POD}}^{HF}$.

The entire data set contains $12000$ samples, and for each $12$ samples a high-fidelity model solution is calculated. Only one high-fidelity model simulation is required for computing the reduced solutions errors for the parametric configuration $\mu$ using reduced-order models of various $K_{POD}$ constructed based on a single high-fidelity trajectory described by parameter $\mu_p$. As such, $1000$ high-fidelity simulations were needed to construct the entire data set. High-fidelity simulations are used to accurately calculate the errors associated with the existing reduced-order models for parametric configurations $\mu$.

\begin{figure}[h]
  \centering
  \subfigure[Isocontours of the errors ${{{\varepsilon}}_{\mu,\mu_p,K_{POD}}^{HF}}$]{\includegraphics[scale=0.35]{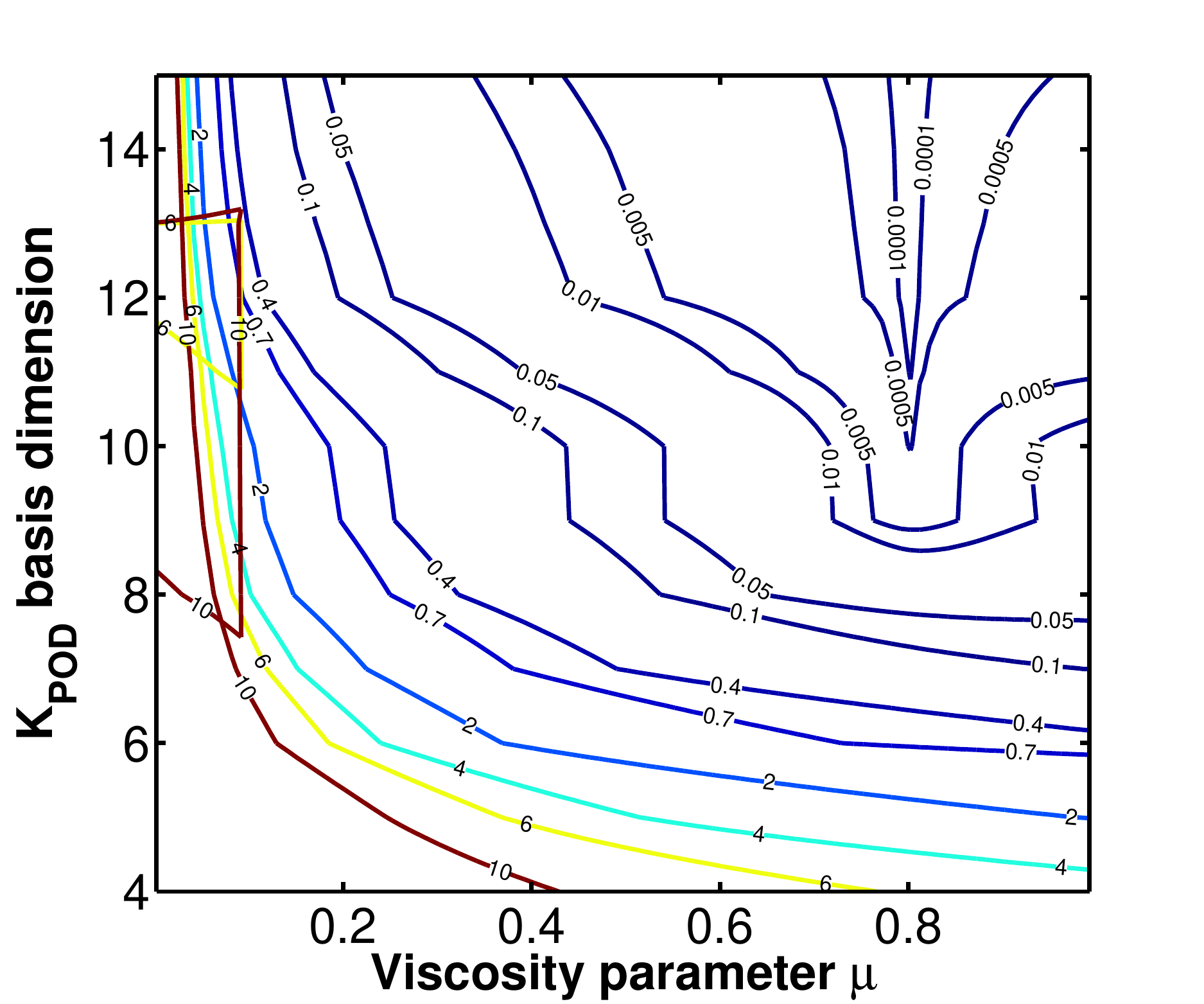}
    \label{fig:parameter_contour_lin}}
    \subfigure[Isocontours for the logarithms of the errors $\log{{{\varepsilon}}_{\mu,\mu_p,K_{POD}}^{HF}}$ ] {\includegraphics[scale=0.35]{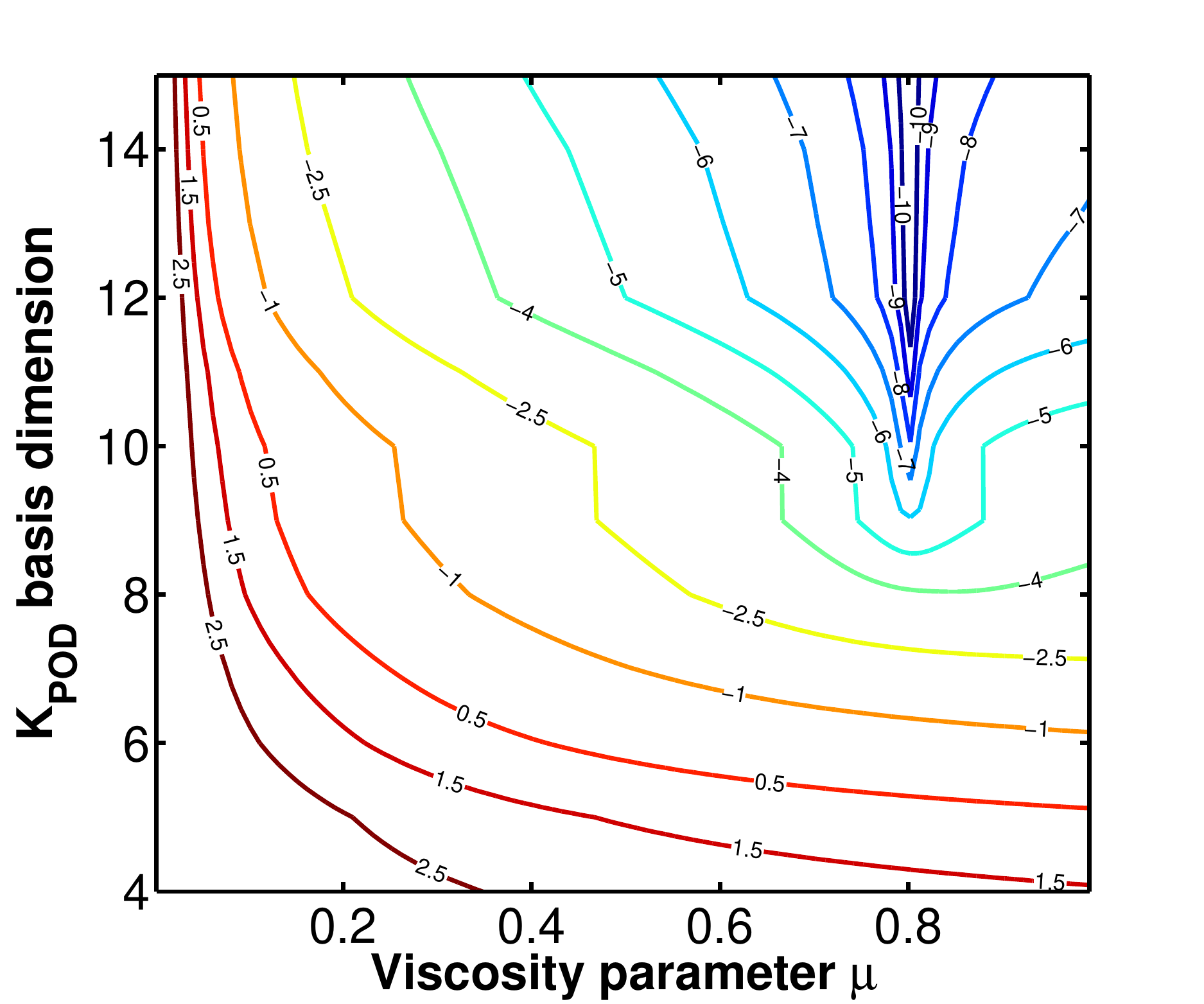}
    \label{fig:parameter_contour_log}}
\caption{Isocontours of the reduced model errors for different POD basis dimensions and parameters $\mu$. The reduced-order model uses a basis constructed from the full order simulation with parameter value $\mu_p=0.8$.}
\label{fig:parameter_contour}
\end{figure}

Figure \ref{fig:parameter_contour} shows isocontours of the error ${{\varepsilon}}_{\mu,\mu_p,K_{POD}}^{HF}$ and $\log{{\varepsilon}}_{\mu,\mu_p,K_{POD}}^{HF}$ of the reduced-order model solution  for various viscosity parameter values $\mu$ and POD basis dimensions. The design of the reduced-order models relies on the high-fidelity trajectory for $\mu_p=0.8$. The target values ${{\varepsilon}}_{\mu,\mu_p,K_{POD}}^{HF}$ vary over a wide range (from  300 to $ 10 ^{-6}$) motivating the choice of implementing models that target  $\log{{{\varepsilon}}_{\mu,\mu_p,K_{POD}}^{HF}}$ to decrease the variance of the predicted results.

%
%


A more detailed analysis, comparing models
 \begin{equation}\label{eqn:prob_model_no_scale}
  \phi_{MP-LROM}^e: \{\mu,\mu_p,K_{POD}\} \mapsto \hat{\varepsilon}_{\mu,\mu_p,K_{POD}}^{HF}
\end{equation}
that target no scaled data and model \eqref{eqn:prob_model_scale} is given in the following.

The approximated MP-LROM models for estimating the local parametric reduced-order model errors are constructed using a Gaussian Process with a squared-exponential covariance kernel \eqref{eq_cov} and a neural network with six hidden layers and hyperbolic tangent sigmoid activation function in each layer. Tables \ref{tab:Param_lin} and \ref{tab:Param_log}  show the averages and variances of errors in prediction of MP-LROM models for different sample sizes. Every subset of samples is selected randomly from a shuffled original data set.
The misfit is computed using the same formulas presented in \eqref{eqn:err_fold} to evaluate the prediction errors. Table \ref{tab:Param_lin} shows the prediction errors of \eqref{eqn:prob_model_no_scale} computed via equation \eqref{eqn:err_fold} with ${y} = {{{\varepsilon}}_{\mu,\mu_p,K_{POD}}^{HF}}$ and $\hat{y} = \hat{\varepsilon}_{\mu,\mu_p,K_{POD}}^{HF}$; i.e., no data scaling; the predictions have a large variance and  a low accuracy. Scaling the data and targeting $\log{{{\varepsilon}}_{\mu,\mu_p,K_{POD}}^{HF}}$ results using \eqref{eqn:prob_model_scale}, reduce the variance of the predictions, and increase the accuracy, as shown in Table \ref{tab:Param_log}. The same formula \eqref{eqn:err_fold} with ${y} = {\log{{\varepsilon}}_{\mu,\mu_p,K_{POD}}^{HF}}$ and $\hat{y} = \widehat{\log{{\varepsilon}}}_{\mu,\mu_p,K_{POD}}^{HF}$ was applied. We notice that, for increasing sample sizes less or equal than $700$ and for scaled data, the variances of GP
and ANN predictions are not necessarily decreasing. This behavior changes and the variances of both regression models decrease for increasing sample sizes larger than $700$ as seen in Table \ref{tab:Param_log}.
The performance of the ANN and GP is highly dependent on the number of samples in the data set.
As the number of data points grows, the accuracy increases and the variance decreases.
The results show that GP outperforms ANN for small numbers of samples
$\leq 1000 $ whereas, for larger data sets, ANN is more accurate than GP.

%
\begin{table}[H]
\begin{center}
    \begin{tabular}{ | l | l | l |  l | l |}
    \hline
     & \multicolumn{2}{|c|}{GP MP-LROM} & \multicolumn{2}{|c|}{ANN MP-LROM} \\
 \hline
Sample size &  $\textnormal{E}_{\rm fold}$   &  $\textnormal{VAR}_{\rm fold}$    & $\textnormal{E}_{\rm fold}$   &  $\textnormal{VAR}_{\rm fold}$
  \\ \hline
 100 & $ 13.4519 $ & $ 5.2372 $ & $  12.5189 $ & $  25.0337 $
     \\ \hline
  400 & $ 6.8003 $ & $ 31.0974 $ & $ 6.9210 $ & $ 26.1814 $
     \\ \hline
   700 & $ 5.6273 $ & $ 14.3949 $ & $ 7.2325 $ & $ 19.9312 $
     \\ \hline
 1000 & $ 3.7148 $ & $ 13.8102 $ & $ 5.6067 $ & $ 14.6488 $
 \\ \hline
 3000 & $ 0.5468 $ & $ 0.0030 $ & $ 1.2858 $ & $ 1.2705 $
 \\ \hline
 5000 & $ 6.0563 $ & $ 22.7761 $ & $ 3.8819 $ & $ 23.9059 $
 \\ \hline
     \end{tabular}
\end{center}
  \caption{Average and variance of error in predictions of MP-LROM models \eqref{eqn:prob_model_no_scale} constructed via ANN and GP using errors ${{\varepsilon}}_{\mu,\mu_p,K_{POD}}^{HF}$ in training data for different sample sizes. \label{tab:Param_lin}}
\end{table}
%

%
\begin{table}[H]
\begin{center}
    \begin{tabular}{ | l | l | l |  l | l |}
    \hline
     & \multicolumn{2}{|c|}{GP MP-LROM} & \multicolumn{2}{|c|}{ANN MP-LROM} \\
 \hline
Sample size &  $\textnormal{E}_{\rm fold}$   &  $\textnormal{VAR}_{\rm fold}$    & $\textnormal{E}_{\rm fold}$   &  $\textnormal{VAR}_{\rm fold}$
     \\ \hline
 100 & $ 0.5319  $ & $  0.0118 $ & $  1.2177 $ & $ 0.1834 $
     \\ \hline
  400 & $ 0.3906$ & $ 0.0007 $ & $ 0.8988 $ & $ 0.2593 $
     \\ \hline
   700 & $ 0.3322 $ & $ 0.0018 $ & $ 0.7320 $ & $ 0.5602 $
     \\ \hline
 1000 & $ 0.2693 $ & $ 0.0002 $ & $ 0.5866 $ & $ 0.4084 $
 \\ \hline
 3000 & $ 0.1558 $ & $  0.5535 \times 10^{-4} $ & $ 0.01202 $ & $  0.2744 \times 10^{-4} $
 \\ \hline
 5000 & $ 0.0775 $ & $  0.4085 \times 10^{-5}  $ & $  0.0075 $ & $ 0.3812 \times 10^{-5}  $
 \\ \hline
     \end{tabular}
\end{center}
 \caption{Average and variance of error in predictions of MP-LROM models \eqref{eqn:prob_model_no_scale} constructed via ANN and GP using logarithms of errors $\log{{{\varepsilon}}_{\mu,\mu_p,K_{POD}}^{HF}}$ in training data for different sample sizes. \label{tab:Param_log}}
\end{table}

\begin{figure}[h]
  \centering
  \subfigure[$\log{{\varepsilon}}_{\mu,\mu_p,K_{POD}}^{HF} - \widehat{\log{{\varepsilon}}}_{\mu,\mu_p,K_{POD}}^{HF}$ \newline {}{- 100 samples}]
  {\includegraphics[scale=0.35]{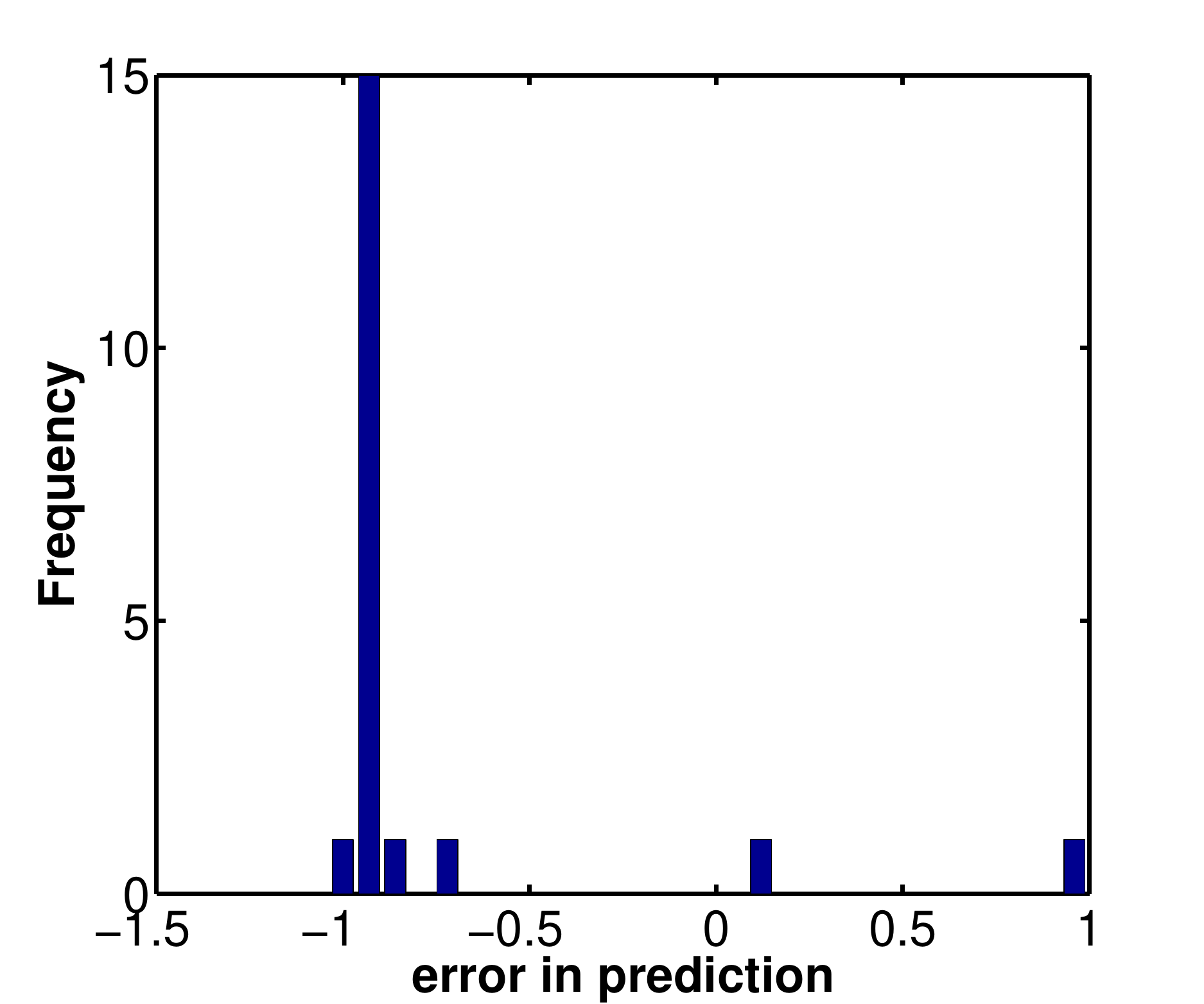}}
  \subfigure[${{\varepsilon}}_{\mu,\mu_p,K_{POD}}^{HF} - {{\hat{\varepsilon}}}_{\mu,\mu_p,K_{POD}}^{HF}$ - 100 samples]{\includegraphics[scale=0.35]
{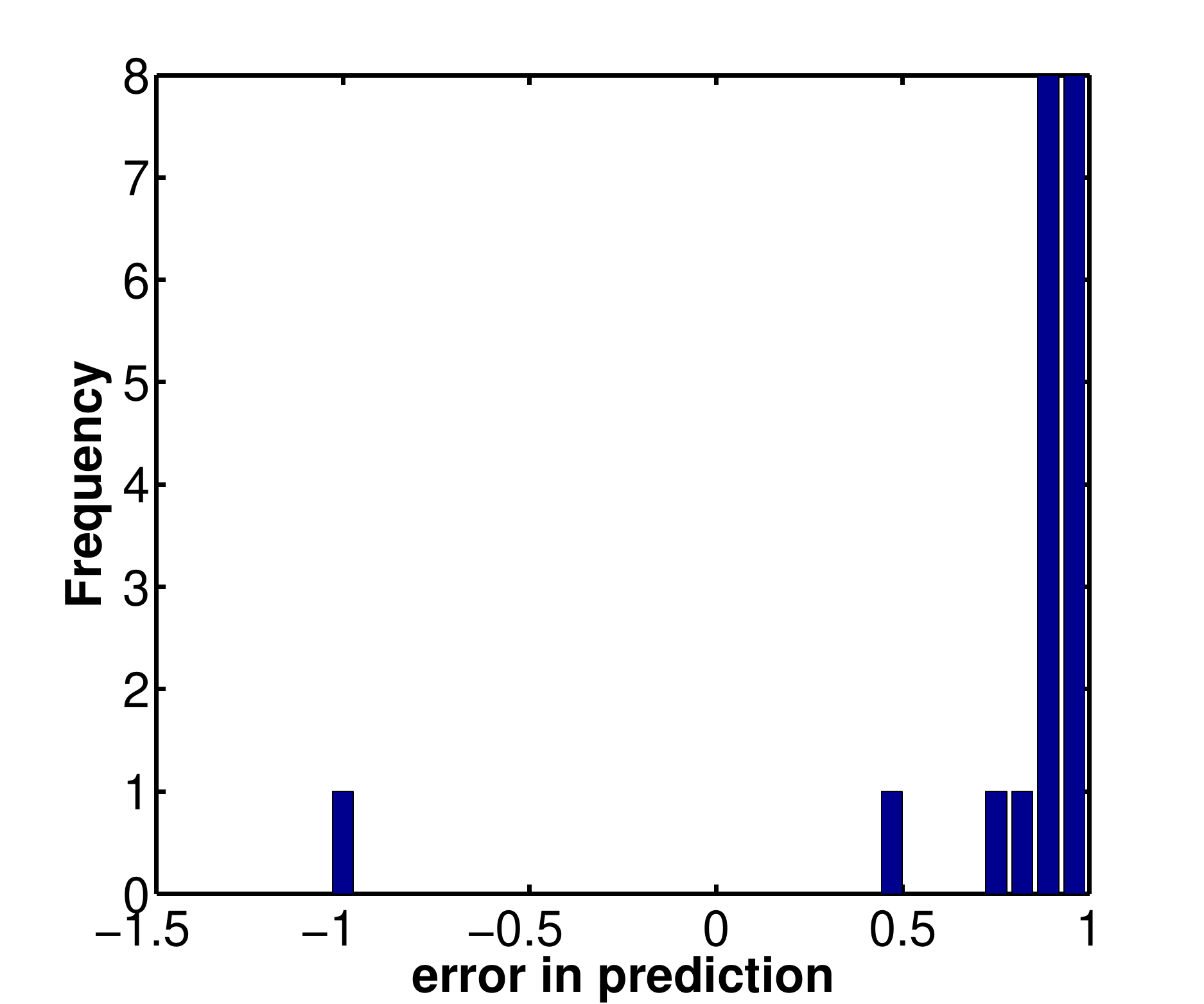}}
  \subfigure[$\log{{\varepsilon}}_{\mu,\mu_p,K_{POD}}^{HF} - \widehat{\log{{\varepsilon}}}_{\mu,\mu_p,K_{POD}}^{HF}$ - 1000 samples] {\includegraphics[scale=0.35]{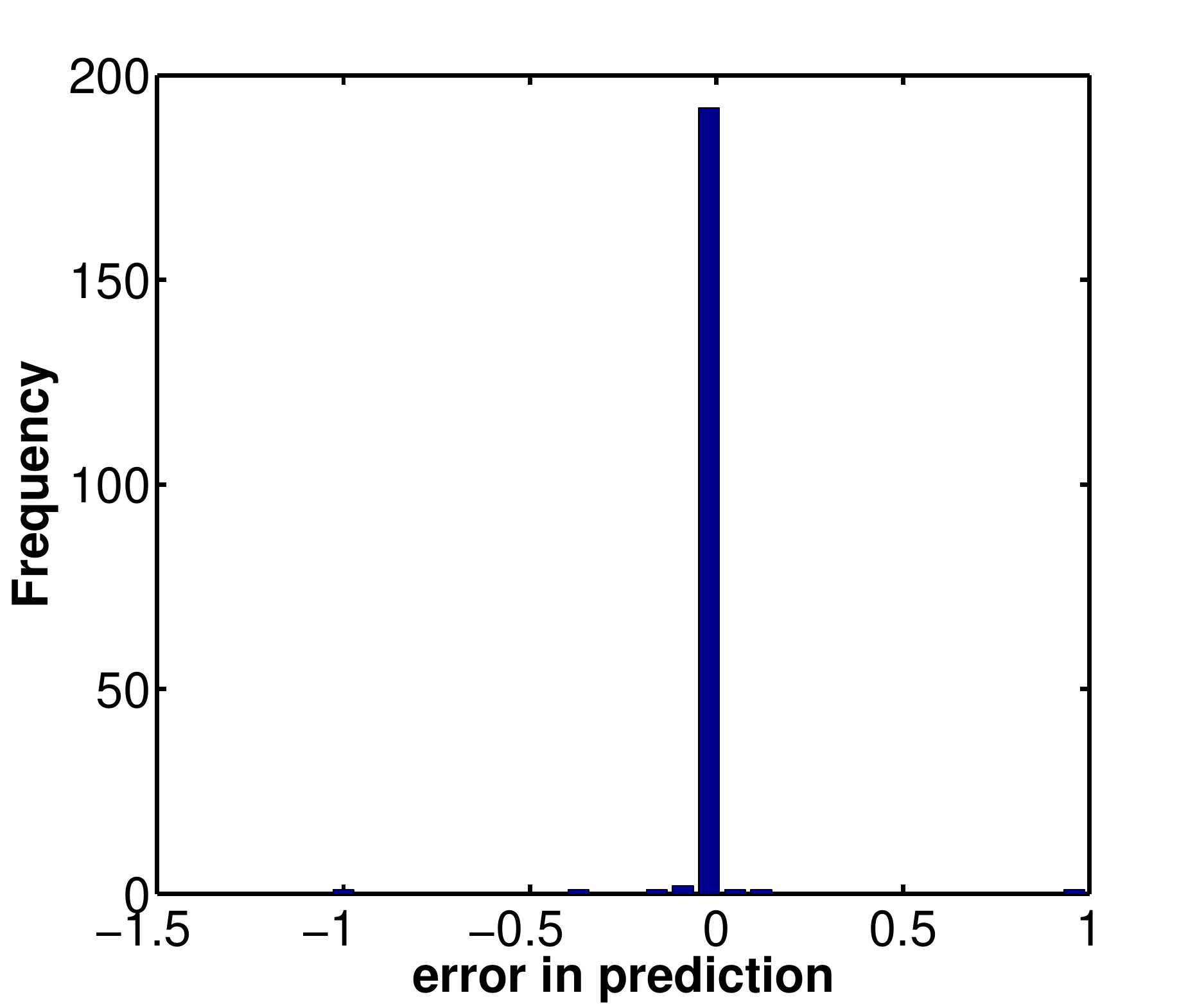}}
  \subfigure[${{\varepsilon}}_{\mu,\mu_p,K_{POD}}^{HF} - {\hat{\varepsilon}}_{\mu,\mu_p,K_{POD}}^{HF}$ - 1000 samples]{\includegraphics[scale=0.35]
{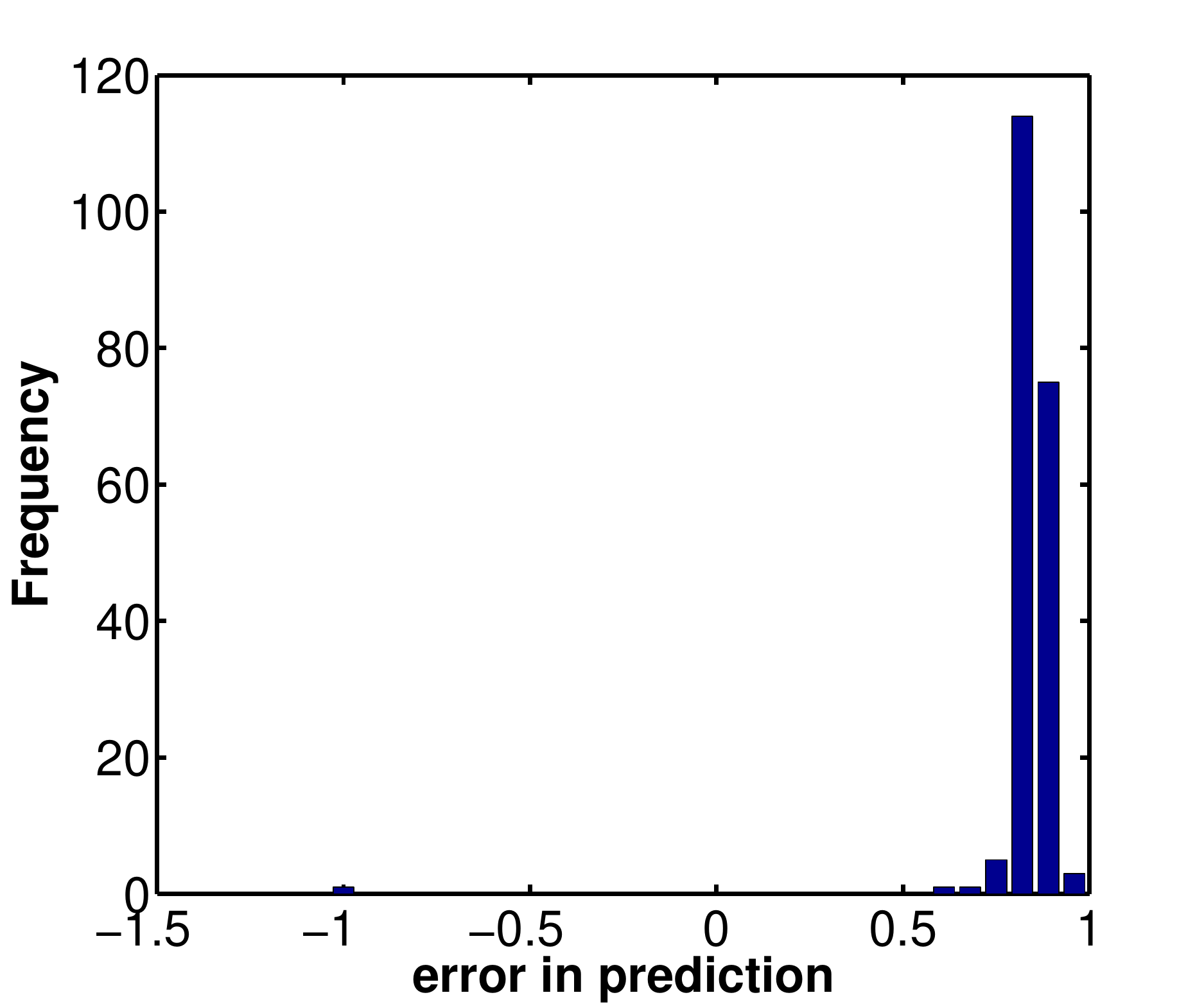}}
\caption{Histogram of errors in prediction using ANN MP-LROM.
\label{fig:ParamHist_NN}}
\end{figure}

\begin{figure}[h]
  \centering
  \subfigure[$\log{{\varepsilon}}_{\mu,\mu_p,K_{POD}}^{HF} - \widehat{\log{{\varepsilon}}}_{\mu,\mu_p,K_{POD}}^{HF}$ - 100 samples] {\includegraphics[scale=0.35]{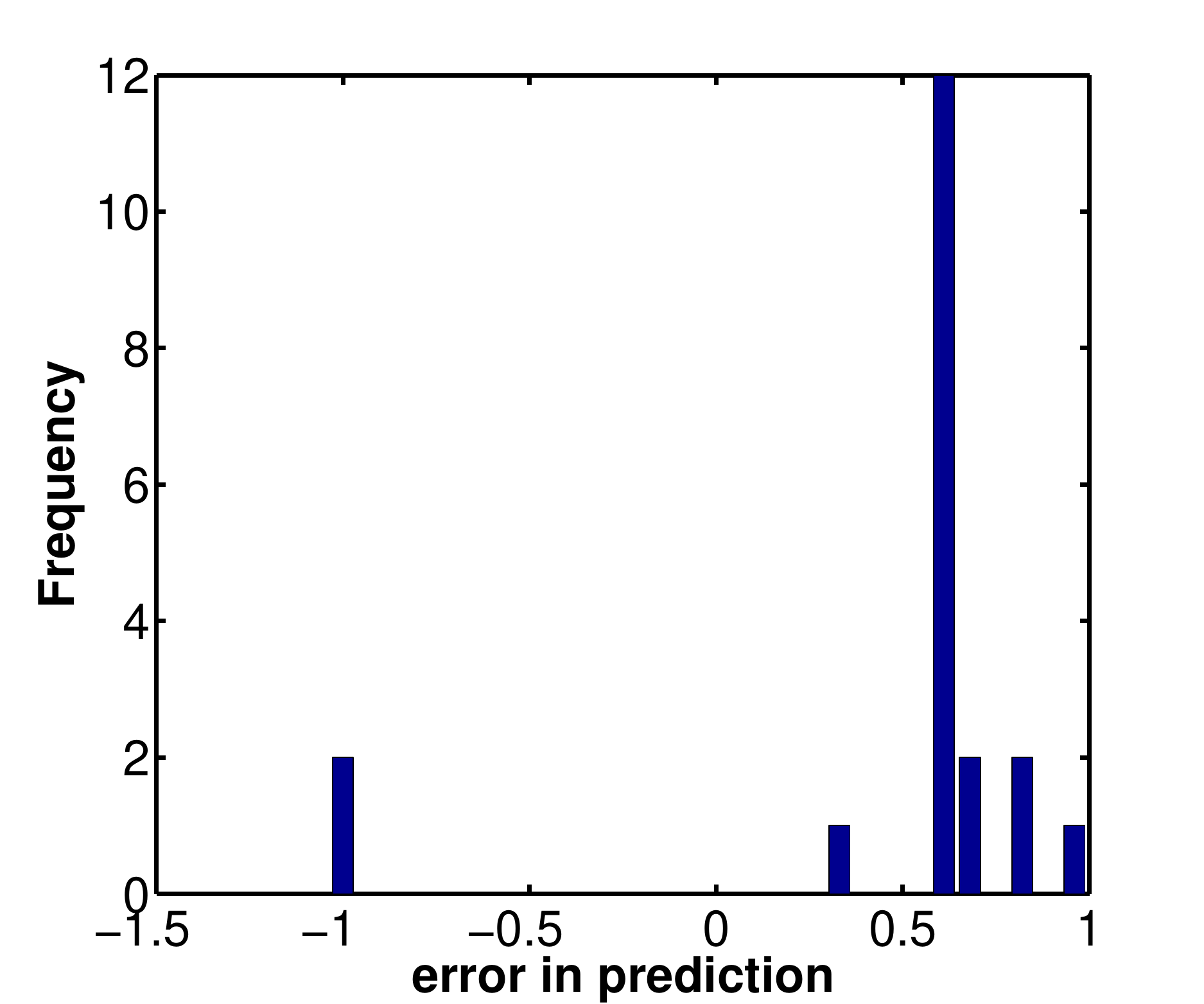}}
  \subfigure[${{\varepsilon}}_{\mu,\mu_p,K_{POD}}^{HF} - {\hat{\varepsilon}}_{\mu,\mu_p,K_{POD}}^{HF}$ - 100 samples]{\includegraphics[scale=0.35]
{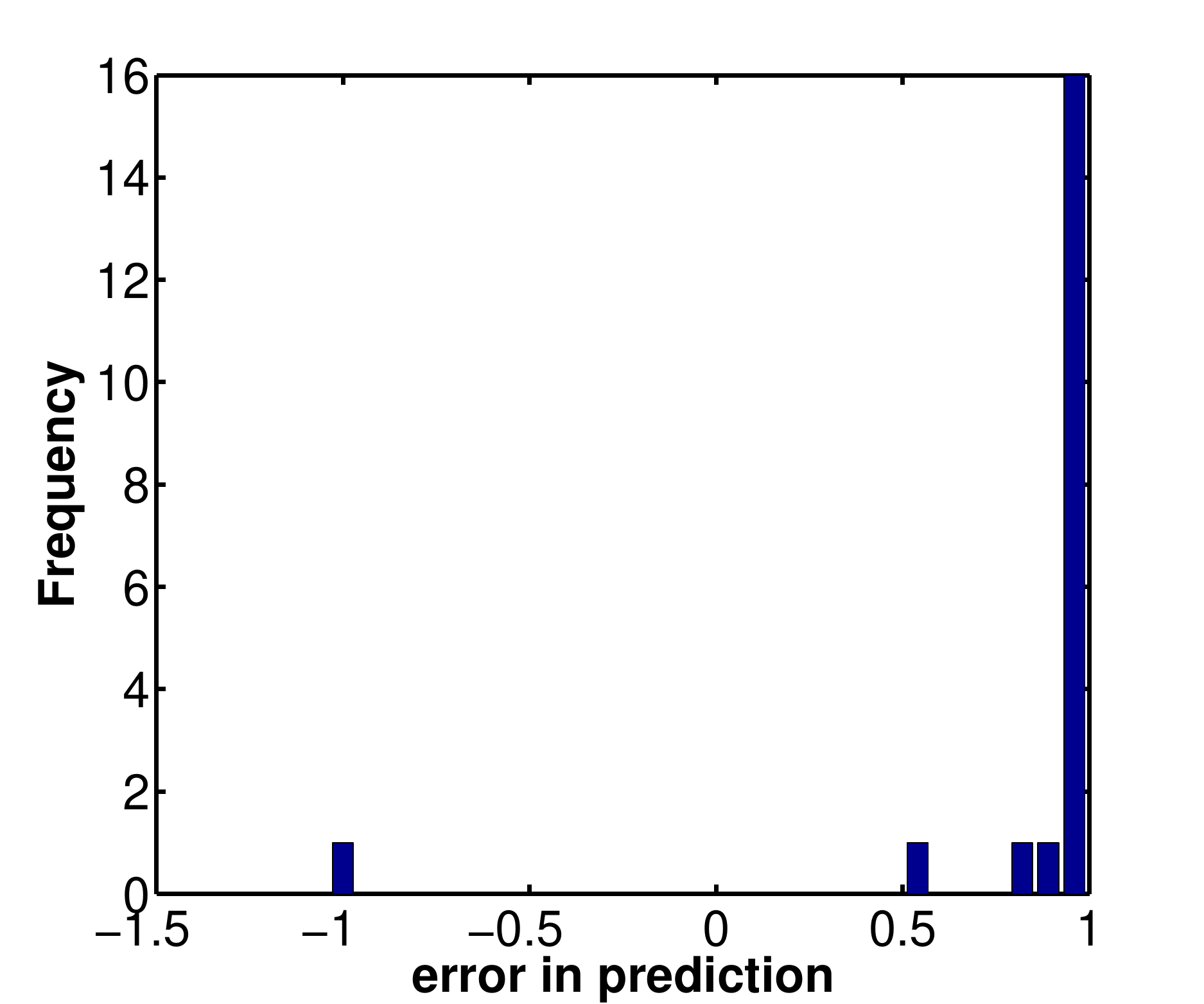}}
  \subfigure[$\log{{\varepsilon}}_{\mu,\mu_p,K_{POD}}^{HF} - \widehat{\log{{\varepsilon}}}_{\mu,\mu_p,K_{POD}}^{HF}$ - 1000 samples]{\includegraphics[scale=0.35]{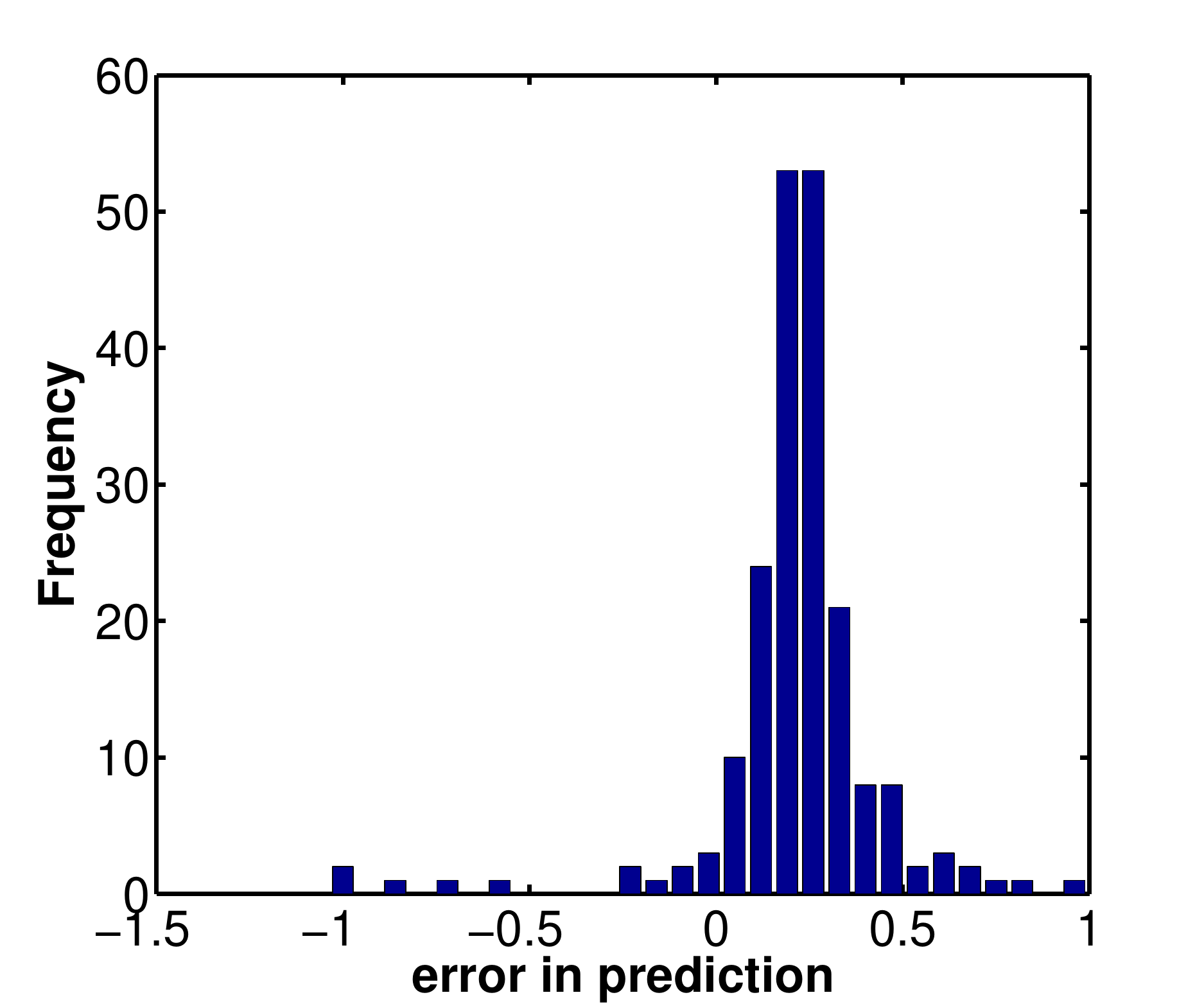}}
  \subfigure[${{\varepsilon}}_{\mu,\mu_p,K_{POD}}^{HF} - {\hat{\varepsilon}}_{\mu,\mu_p,K_{POD}}^{HF}$ - 1000 samples]{\includegraphics[scale=0.35]
  {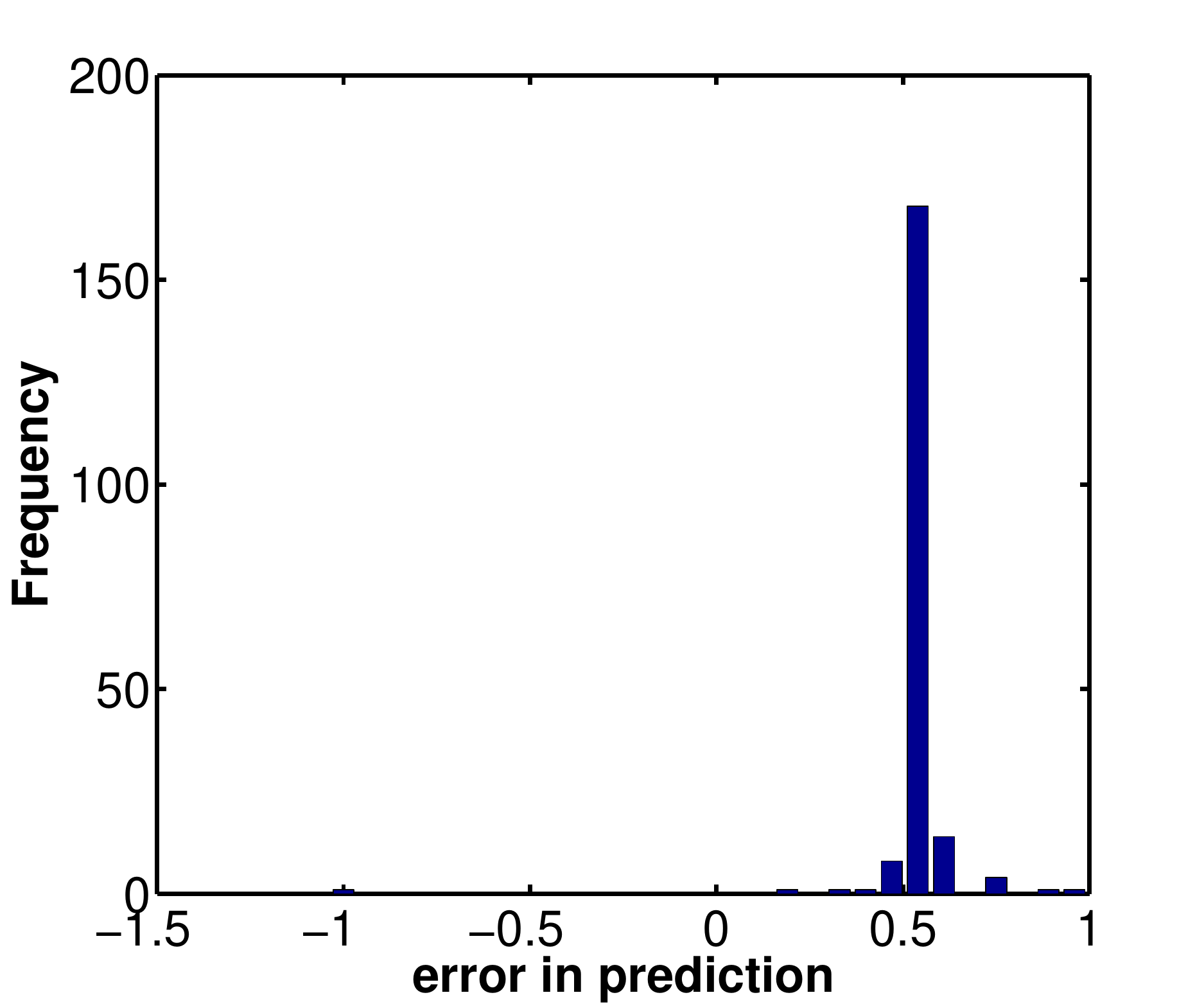}}
\caption{Histogram of errors in prediction using GP MP-LROM.
\label{fig:ParamHist_GP}}
\end{figure}

 Figures \ref{fig:ParamHist_NN} and \ref{fig:ParamHist_GP} show the corresponding histogram of the errors in prediction
 of MP-LROM models \eqref{eqn:prob_model_scale} and \eqref{eqn:prob_model_no_scale}  using $100$ and $1000$ training samples for
 ANN and GP methods, respectively. 
 The histograms shown in Figure \ref{fig:ParamHist_GP} can assess the validity of GP assumptions \eqref{GP_Dist}, \eqref{GP_training}, \eqref{GP_prior}. The difference between the true and estimated values should behave as samples from the distribution $ \mathcal{N} (0, \sigma_n^2) $ \cite{drohmann2015romes}. In our case they are hardly normally distributed and this indicates that the data sets are not from Gaussian distributions.
%
%
%

 Scaling the data and targeting $\log{{{\varepsilon}}_{\mu,\mu_p,K_{POD}}^{HF}}$ errors clearly improve the performance of the MP-LROM models. Consequently for the rest of the manuscript we will only use model \eqref{eqn:prob_model_scale}.

 To asses the quality of the MP-LROM models, we also implemented a five-fold cross-validation test over the entire dataset. The results computed using formula \eqref{eqn:err_fold_average} are shown in Table \ref{tab:experm2}. ANN outperforms GP and estimates the errors more accurately. It also has less variance than the Gaussian Process which indicates it has more stable predictions.

%
\begin{table}[H]
\begin{center}
    \begin{tabular}{ | l | l | l |}
    \hline
   & $\textnormal{E} $  & $ \textnormal{VAR} $
      \\ \hline
ANN MP-LROM  & $0.004004$ & $2.16 \times 10^{-6	}$
     \\ \hline
GP MP-LROM & $0.092352$ & $ 1.32 \times 10^{-5} $
 \\ \hline
     \end{tabular}
\end{center}
 \caption{MP-LROM statistical results  over five-fold cross-validation.}
  \label{tab:experm2}
\end{table}

\begin{figure}[t!]
	\begin{centering}
	\includegraphics[width=0.5\textwidth, height=0.40\textwidth]{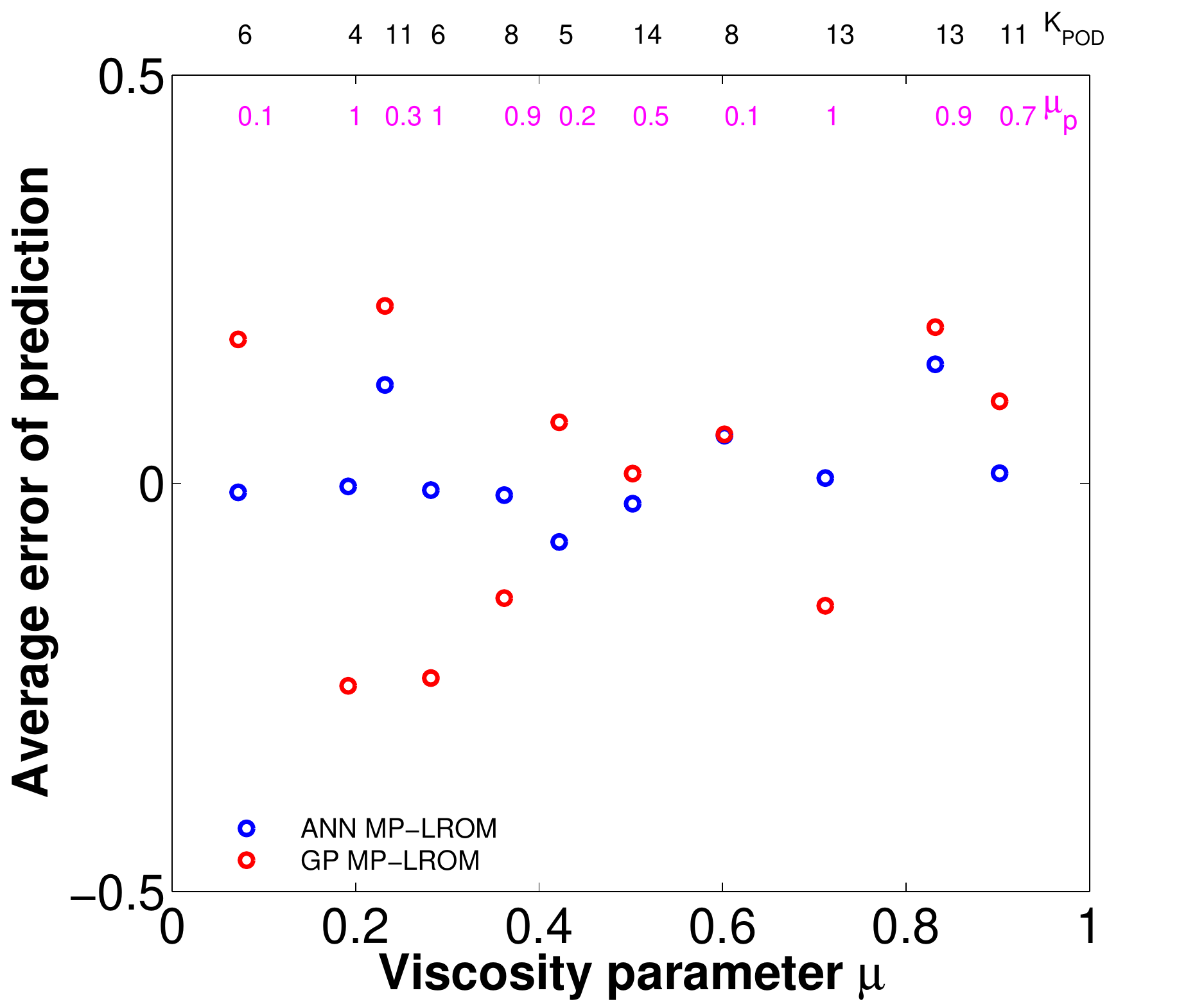}
    \caption{The average of errors in predictions using five different trained models. The top labels shows the corresponding $\mu_p$ and $K_{POD}$ for each parameter $\mu$. }
	\label{fig:expm2_error_estimates}
	\end{centering}
\end{figure}

Figure \ref{fig:expm2_error_estimates} illustrates the average of errors in prediction of five different errors models computed using ANN and GP regression methods. The error models were constructed using a training set formed by $80\%$ randomly selected data of the entire data set. The predictions were made using a fixed test set randomly selected from the entire data set and contains various values of $\mu$, $K_{POD}$ and $\mu_p$ shown in the x-axes of Figure \ref{fig:expm2_error_estimates}.  Building different GP and ANN MP-LROM error models, each trained on different part of the data set and then testing them with the same fixed test set, reduces the bias in prediction. Again, ANN outperforms GP having more accurate errors estimates.

%
%
We also compared the MP-LROM models with those obtained by implementing ROMES method \cite{drohmann2015romes} and MFC technique \cite{alexandrov2001approximation}. The ROMES method constructs univariate models
\begin{equation} \label{eqn::ROMES_math_framework}
  \phi_{ROMES}: \log \rho(\mu) \mapsto \log\varepsilon_{\mu}^{HF},
\end{equation}
where the input $\rho(\mu)$ consists of error indicators. Examples of indicators include residual norms, dual-weighted residuals and other error bounds. MFC implements input-output models
\begin{equation}\label{eqn::MFC_math_framework}
  \phi_{MFC}: \mu \mapsto \log\varepsilon_{\mu}^{HF},
\end{equation}
where the input of error models is the viscosity parameter $\mu$. Both ROMES and MFC methods use a global reduced-order model with a fixed dimension in contrast to our method that employs local reduced-order models with various dimensions. ROMES and MFC models are univariate whereas the MP-LROM models are multivariate.

To accommodate our data set to the requirements of the ROMES and  MFC methods, we separated the data set into multiple subsets. Each of these subsets has $100$ samples corresponding to a single $\mu_p$ and $K_{POD}$ and $100$ values of parameter $\mu \in \{ 0.01,0.02 \ldots, 1\}$, so $100$ high-fidelity simulations are required. For each subset we constructed ANN and GP models to approximate the input-output models defined in \eqref{eqn::ROMES_math_framework} and \eqref{eqn::MFC_math_framework} using the same training set. In the case of ROMES method we employed the logarithms of residuals norms as inputs.  We first computed the corresponding reduced-order solution and then the associated logarithm of residual norm by using the projected  reduced order solution into the high-fidelity model for parameter $\mu$. The output of both ROMES and MFC models approximates the logarithm of the Frobenius norm of the reduced-order-model errors.

\begin{figure}[h]
  \centering
  \subfigure[MFC] {\includegraphics[scale=0.24]{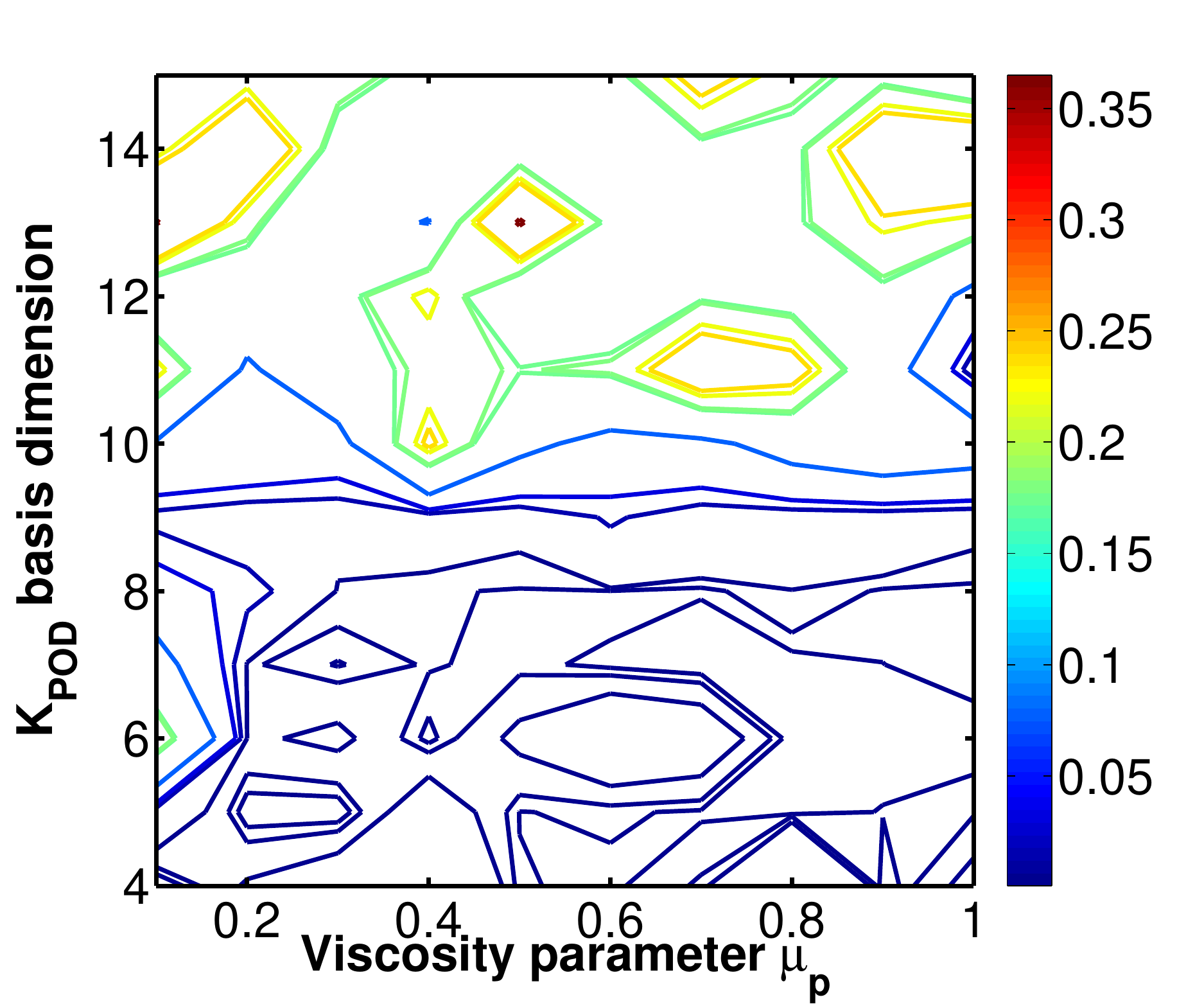}}
  \subfigure[ROMES ]{\includegraphics[scale=0.24]{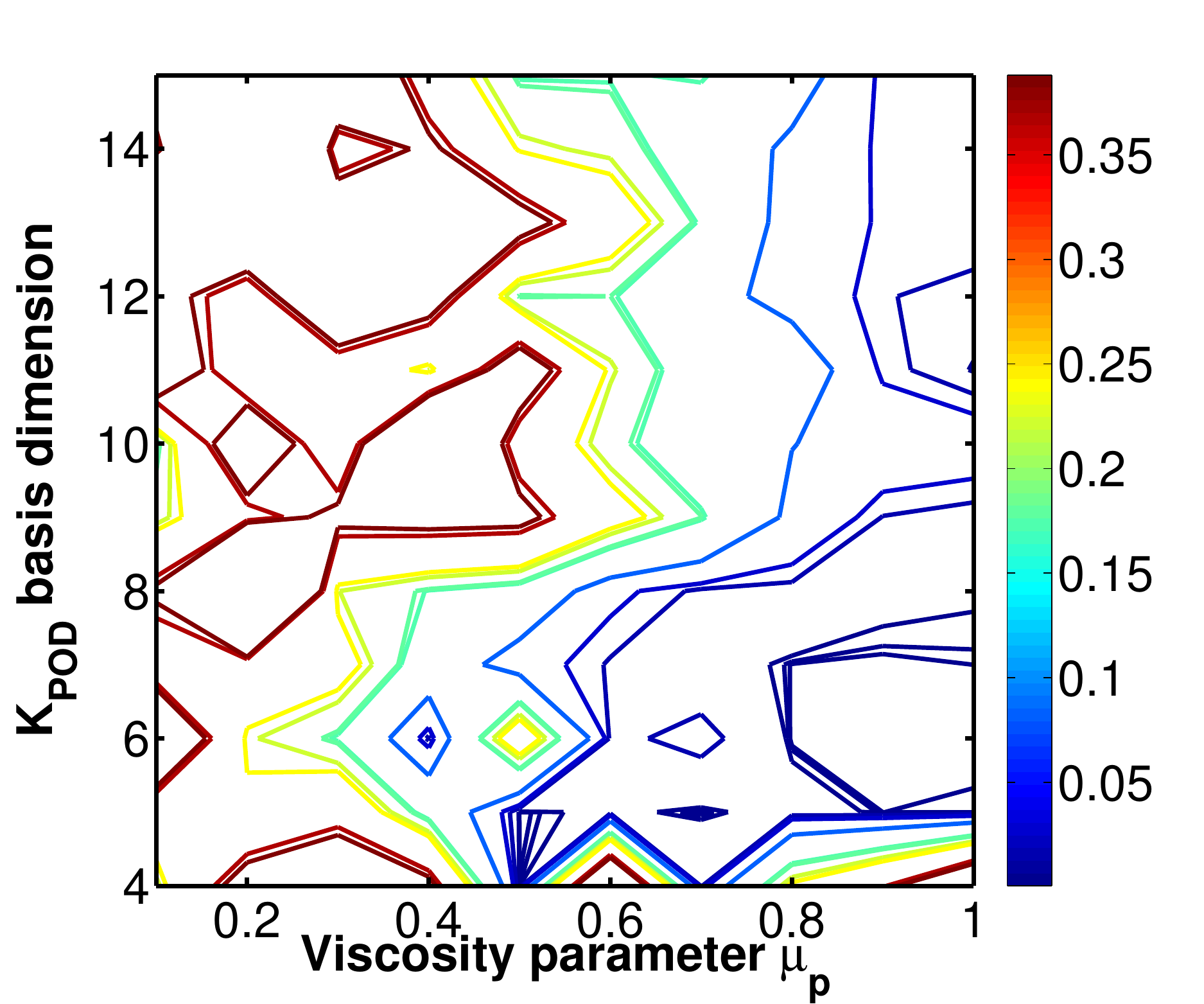}}
  \subfigure[MP-LROM ]{\includegraphics[scale=0.24]{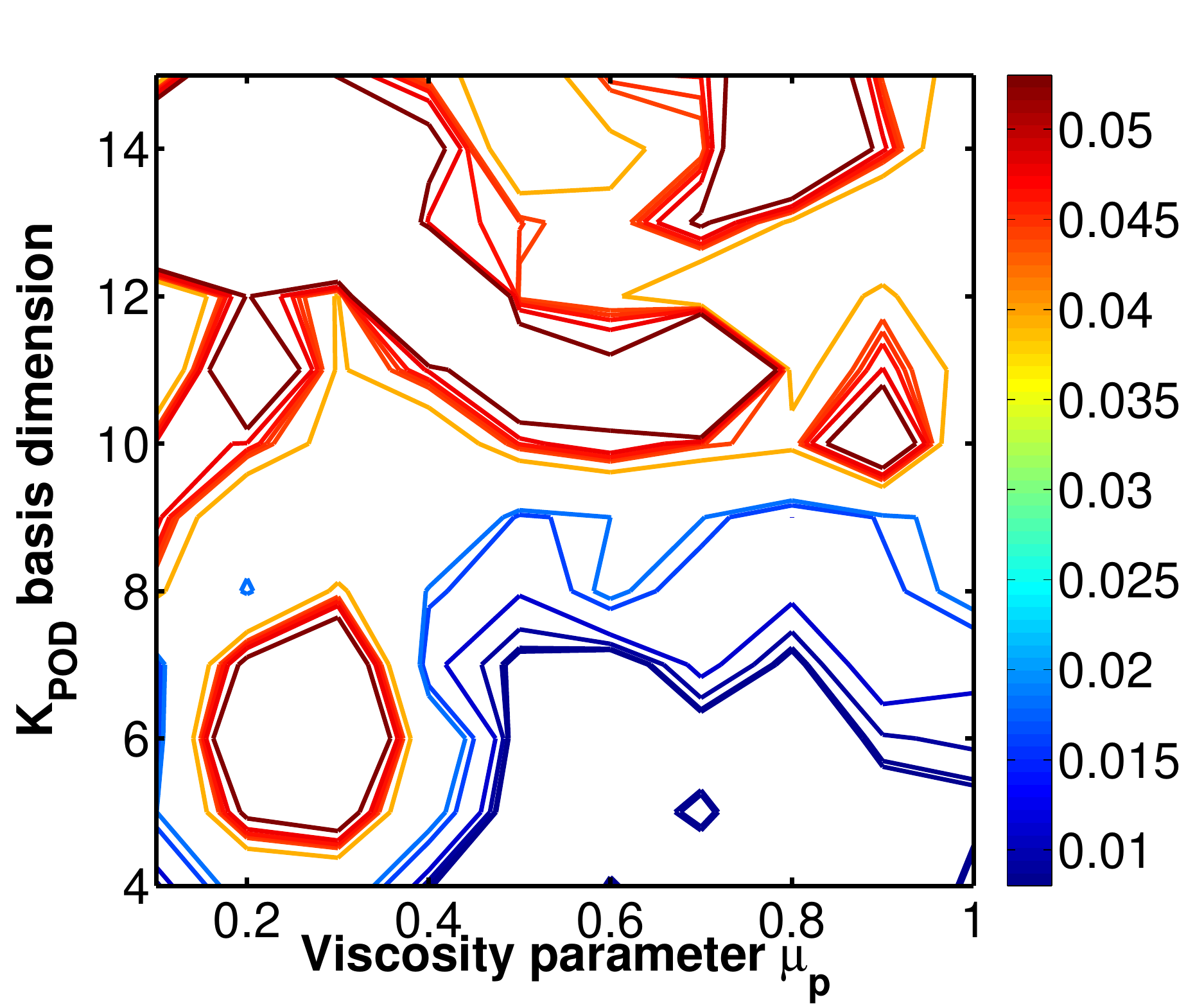}}
\caption{Isocontours for the $\textnormal{E}_{\rm fold}$ using GP method.
\label{fig:contours_GP}}
\end{figure}

\begin{figure}[h]
  \centering
  \subfigure[MFC] {\includegraphics[scale=0.24]{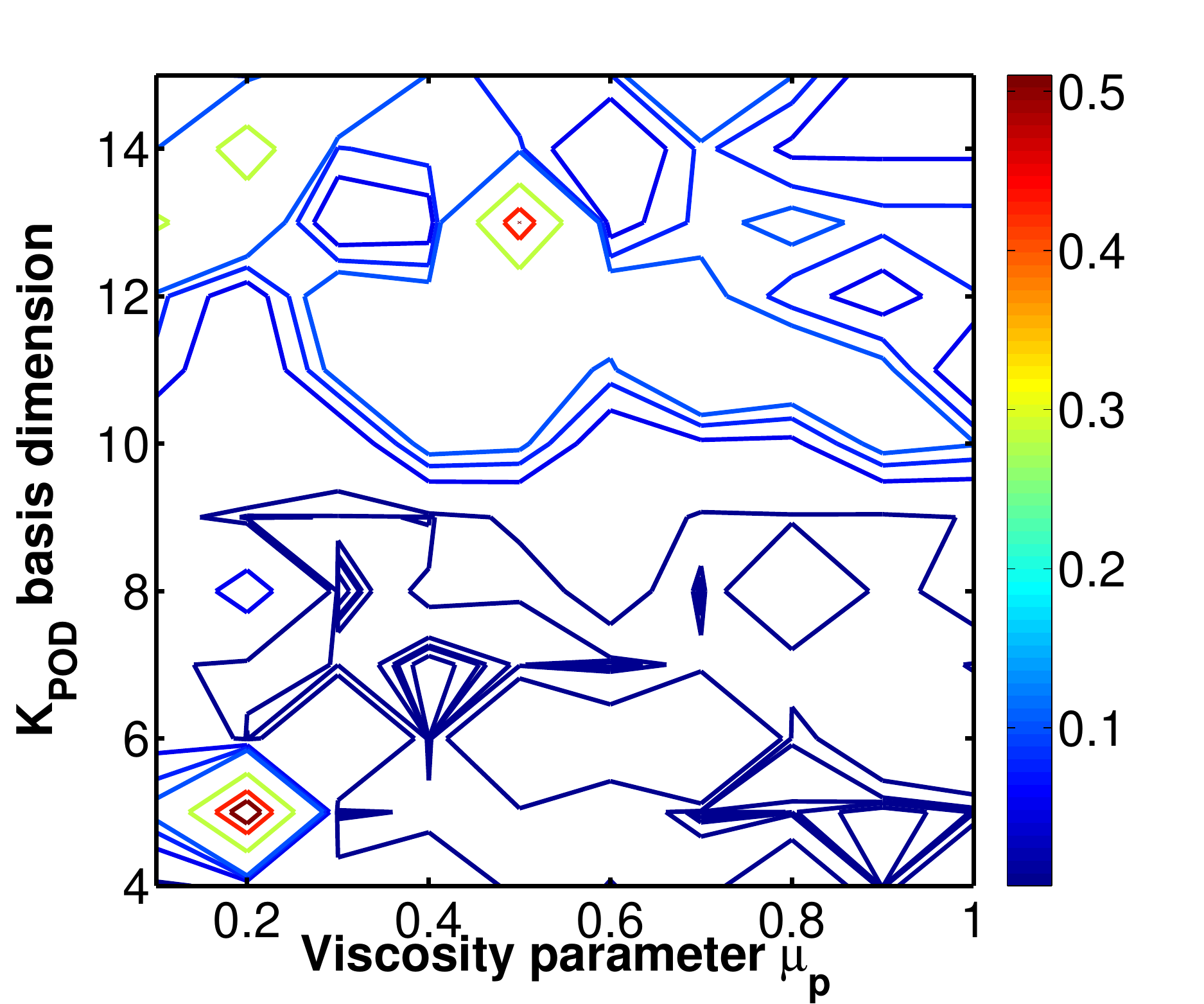}}
  \subfigure[ROMES ]{\includegraphics[scale=0.24]{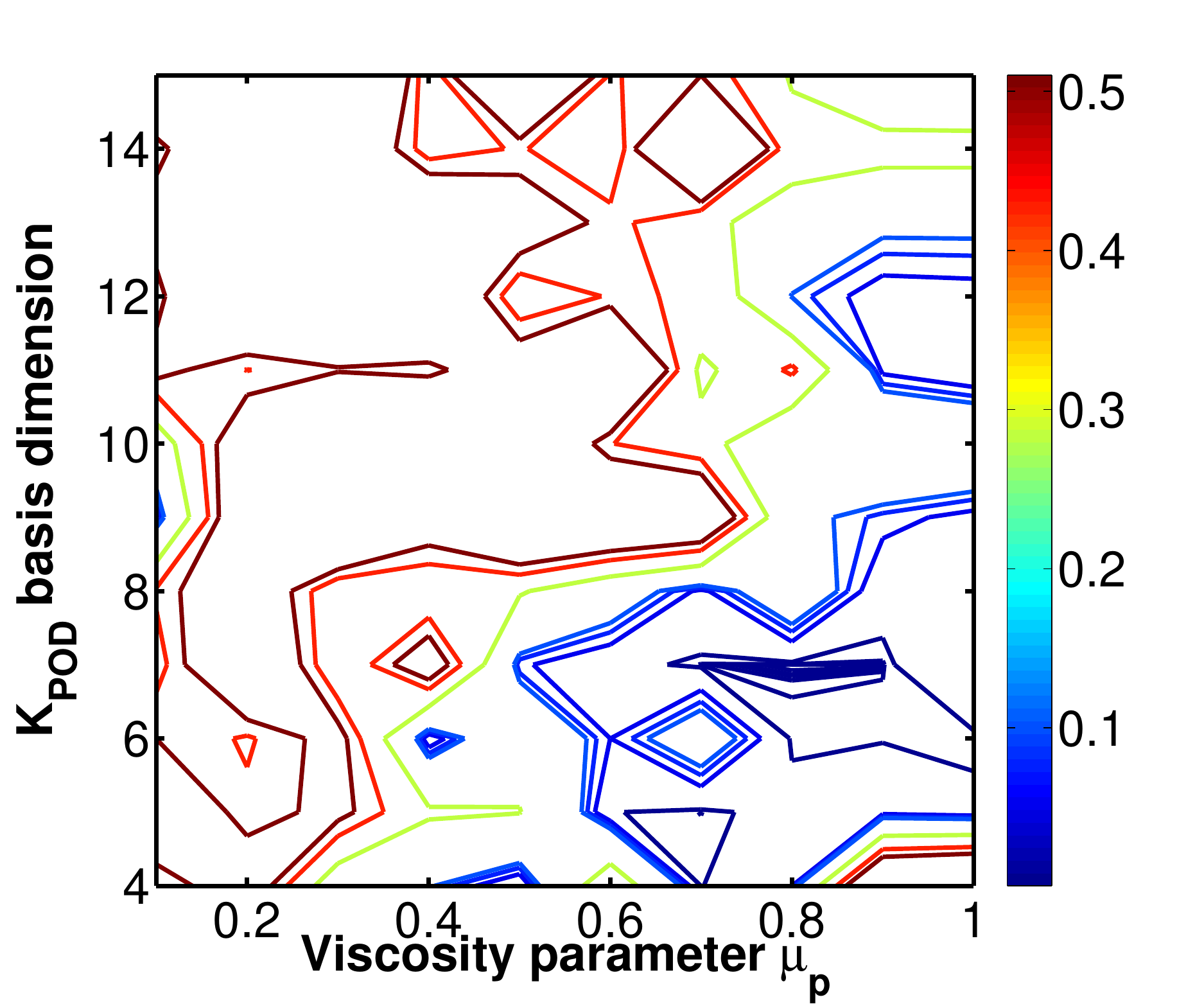}}
   \subfigure[MP-LROM ]{\includegraphics[scale=0.24]{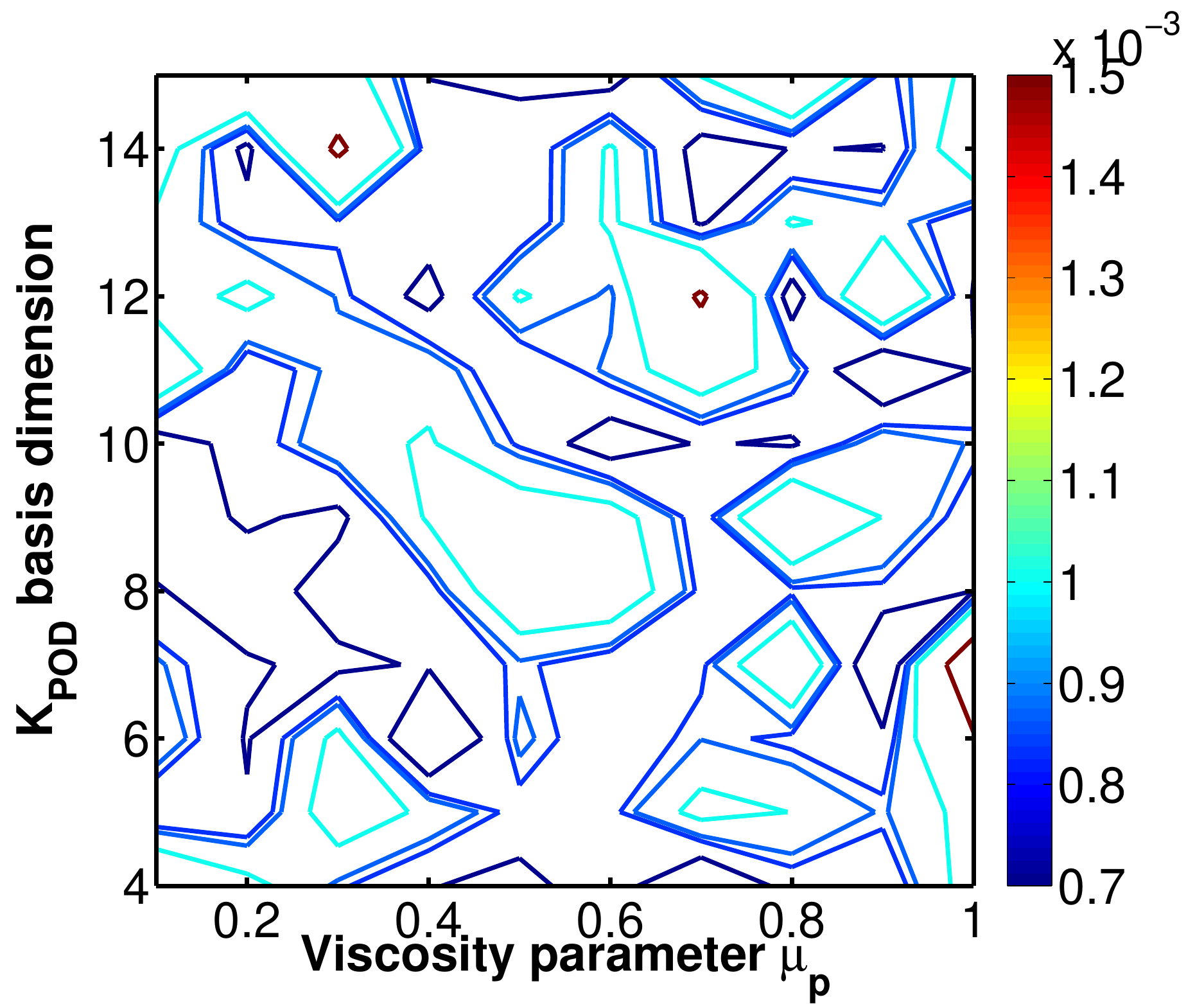}}
\caption{Isocontours for the $\textnormal{E}_{\rm fold}$ using ANN method.
\label{fig:contours_NN}}
\end{figure}

\begin{figure}[h]
  \centering
  \subfigure[MFC] {\includegraphics[scale=0.24]{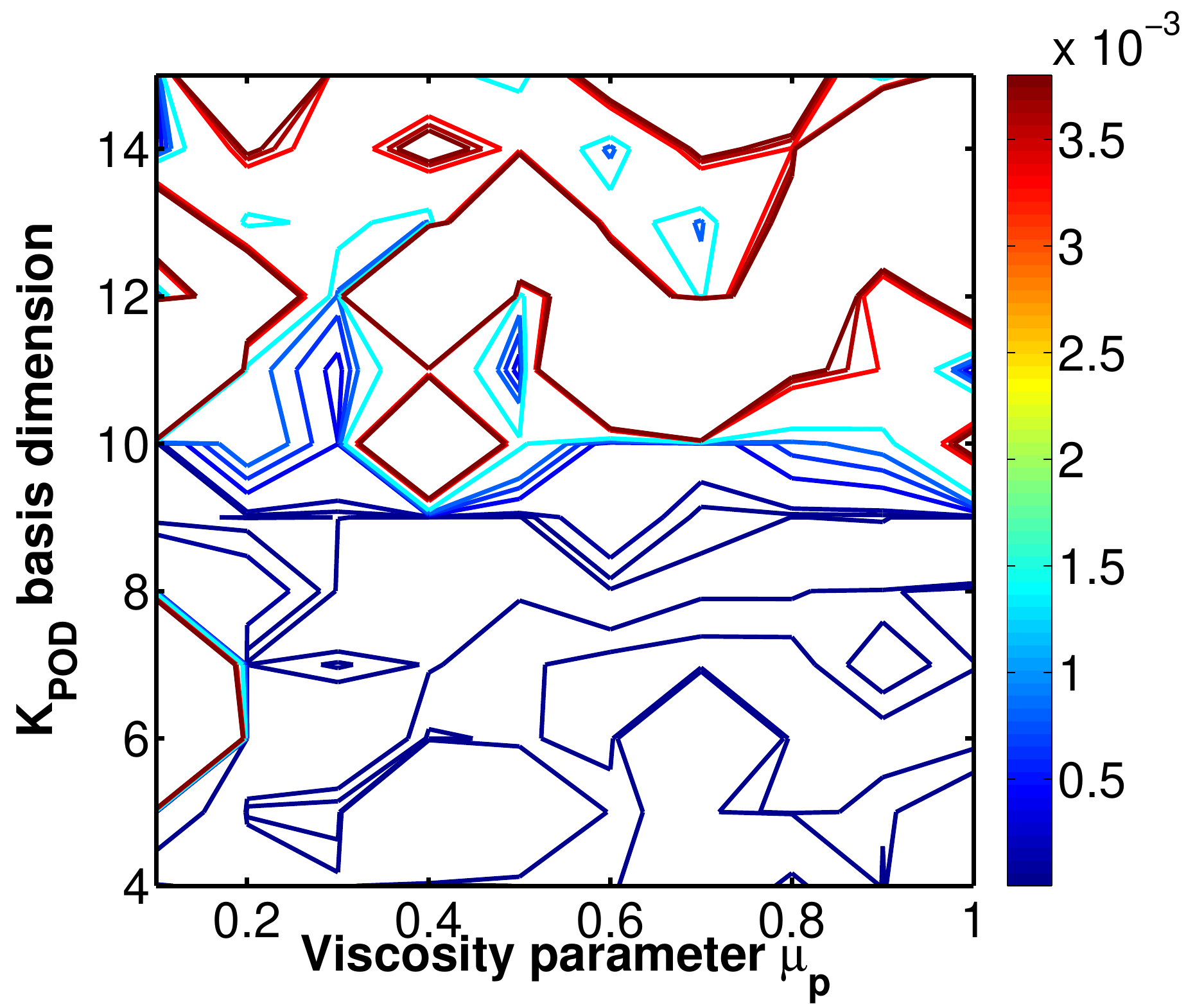}}
  \subfigure[ROMES]{\includegraphics[scale=0.24]{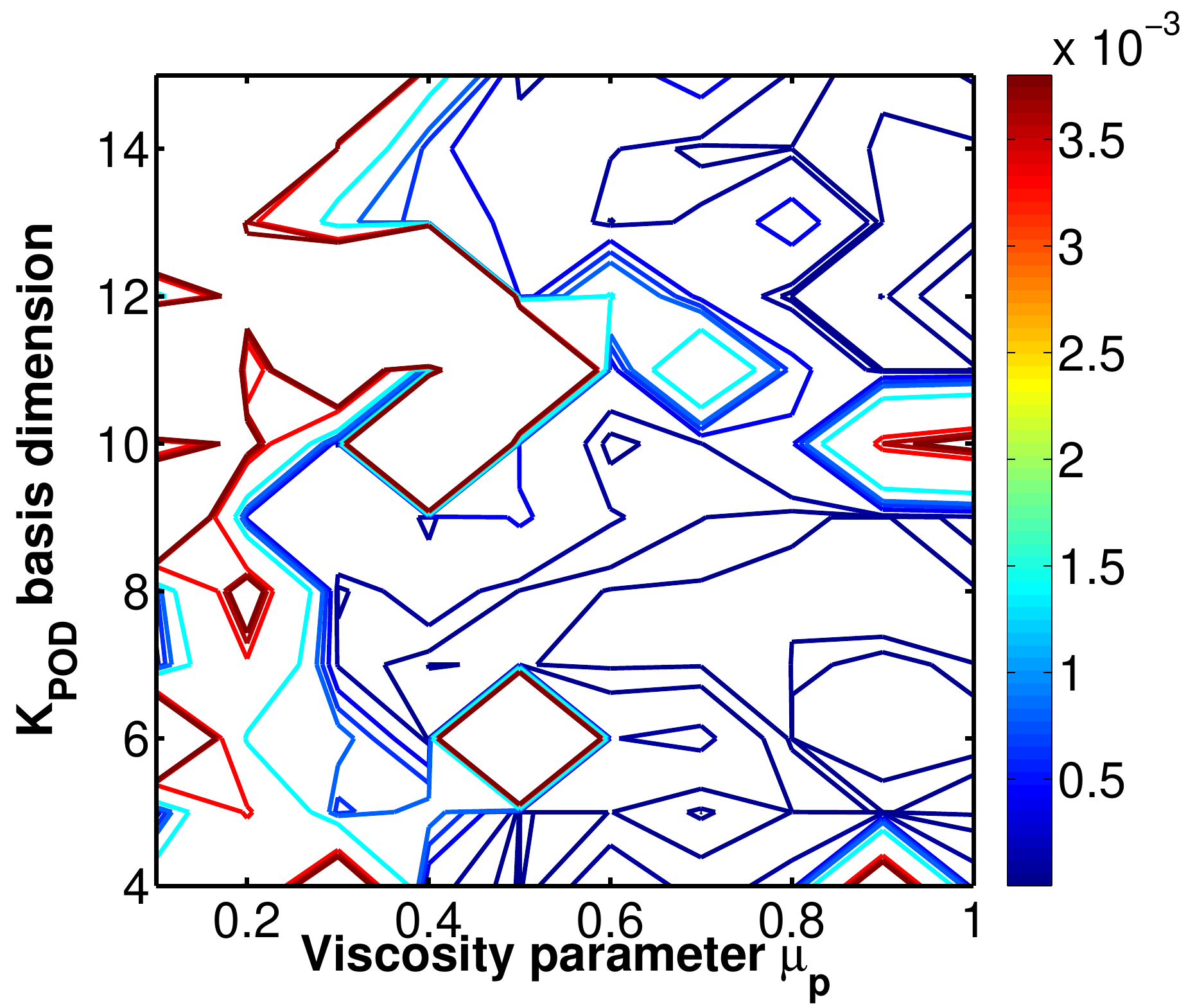}}
  \subfigure[MP-LROM ]{\includegraphics[scale=0.24]{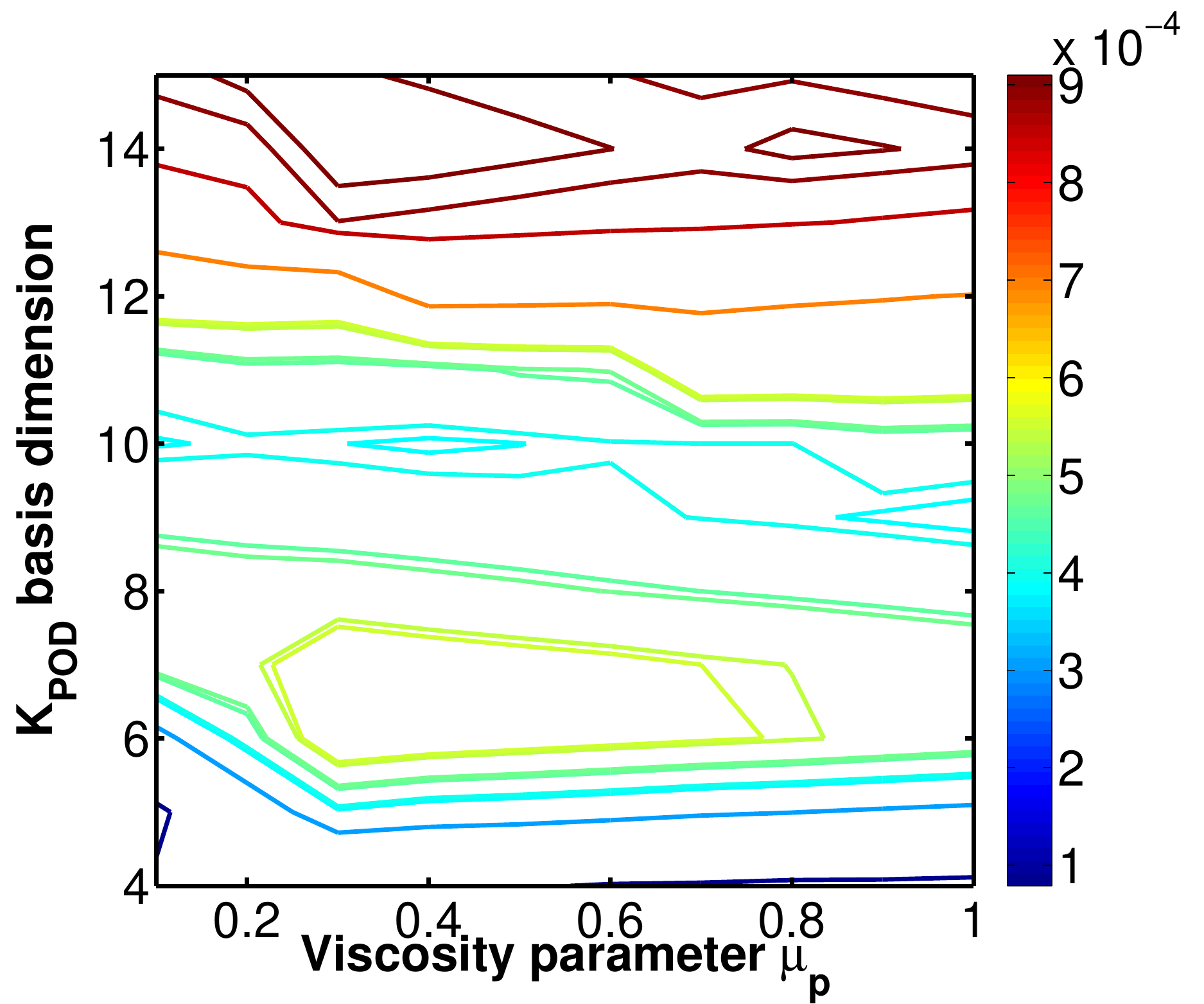}}
\caption{Isocontours for the $\textnormal{VAR}_{\rm fold}$ using GP method.
\label{fig:VARcontours_GP}}
\end{figure}

\begin{figure}[h]
  \centering
  \subfigure[MFC] {\includegraphics[scale=0.24]{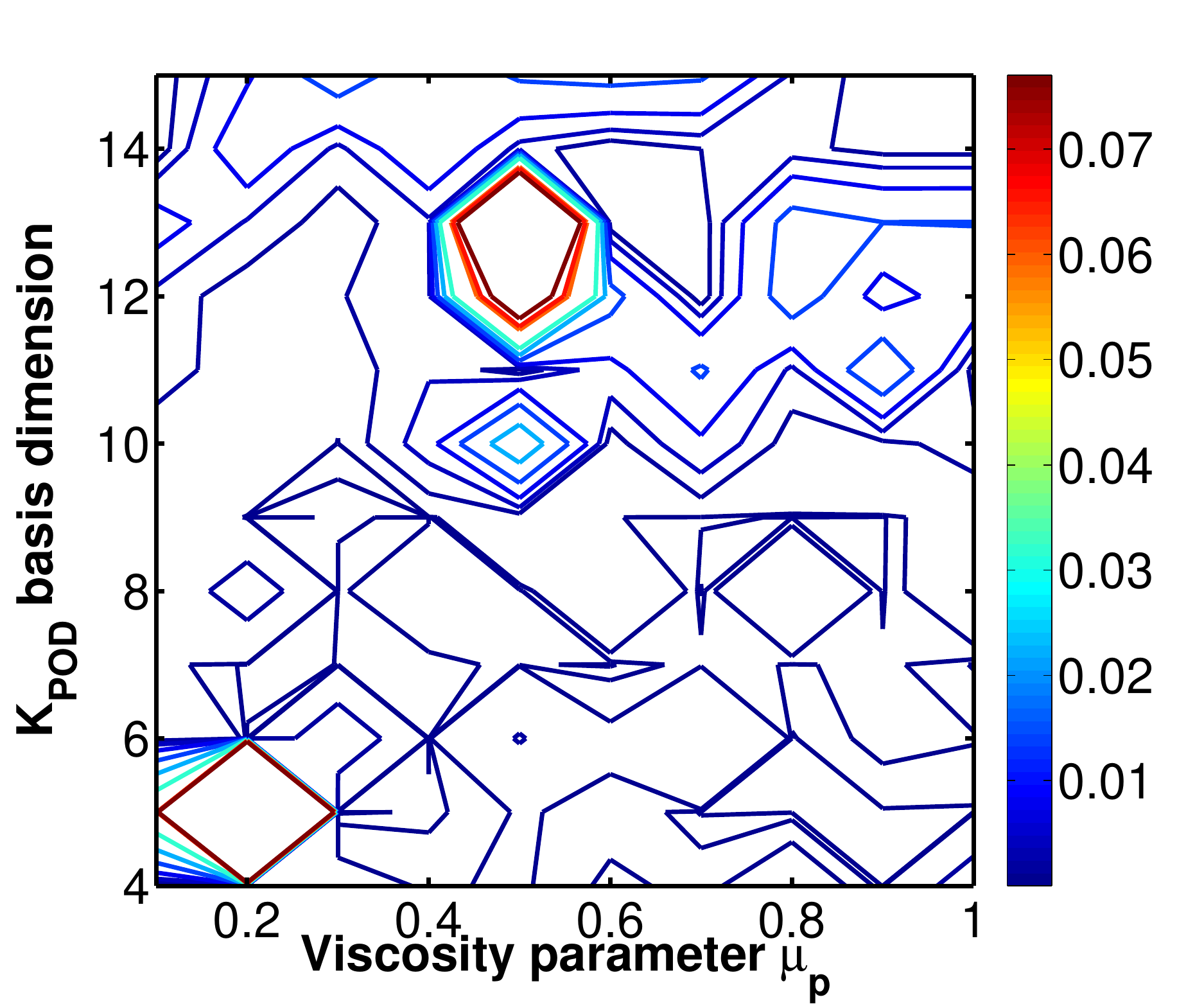}}
  \subfigure[ROMES ]{\includegraphics[scale=0.24]{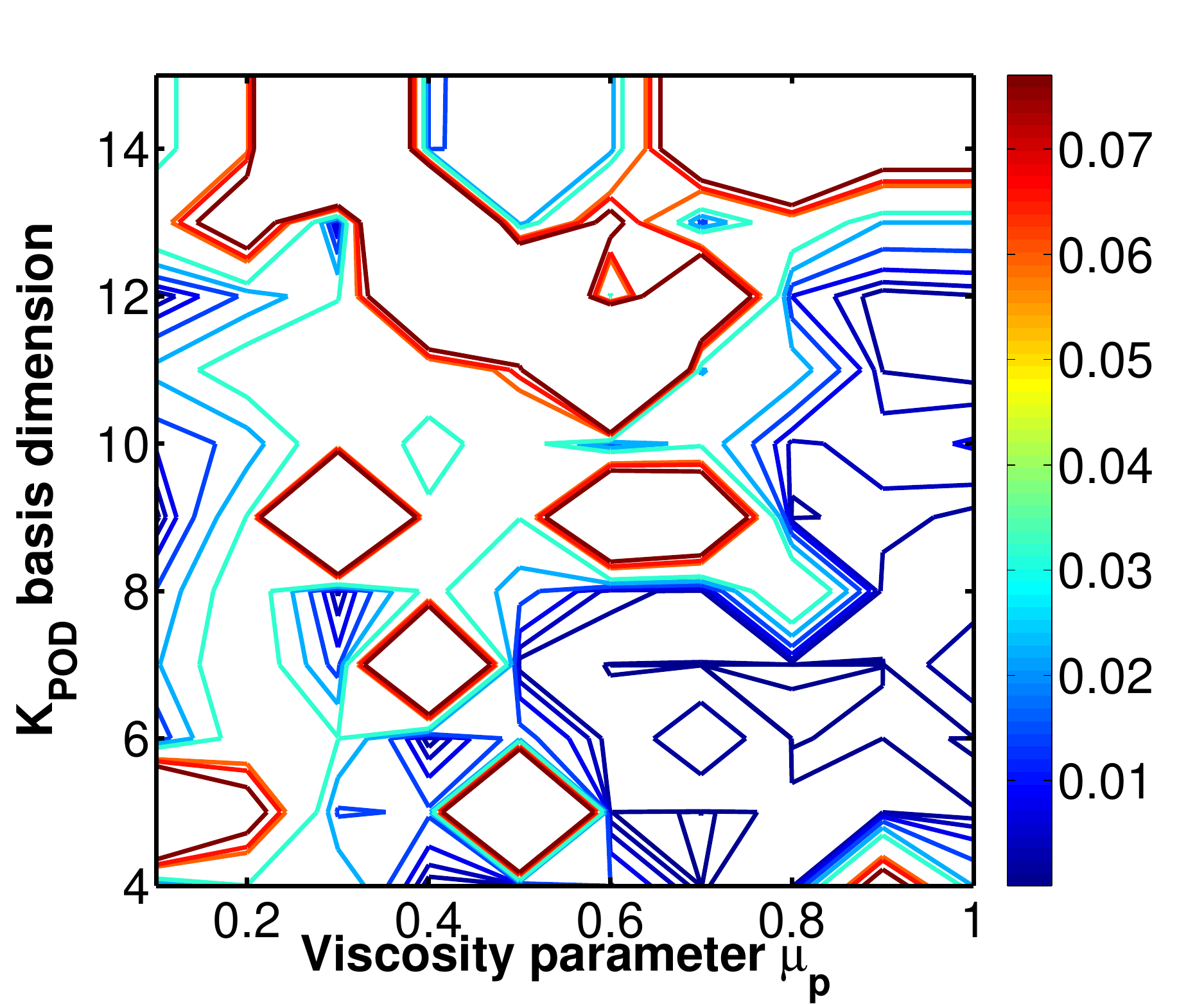}}
  \subfigure[MP-LROM]{\includegraphics[scale=0.24]{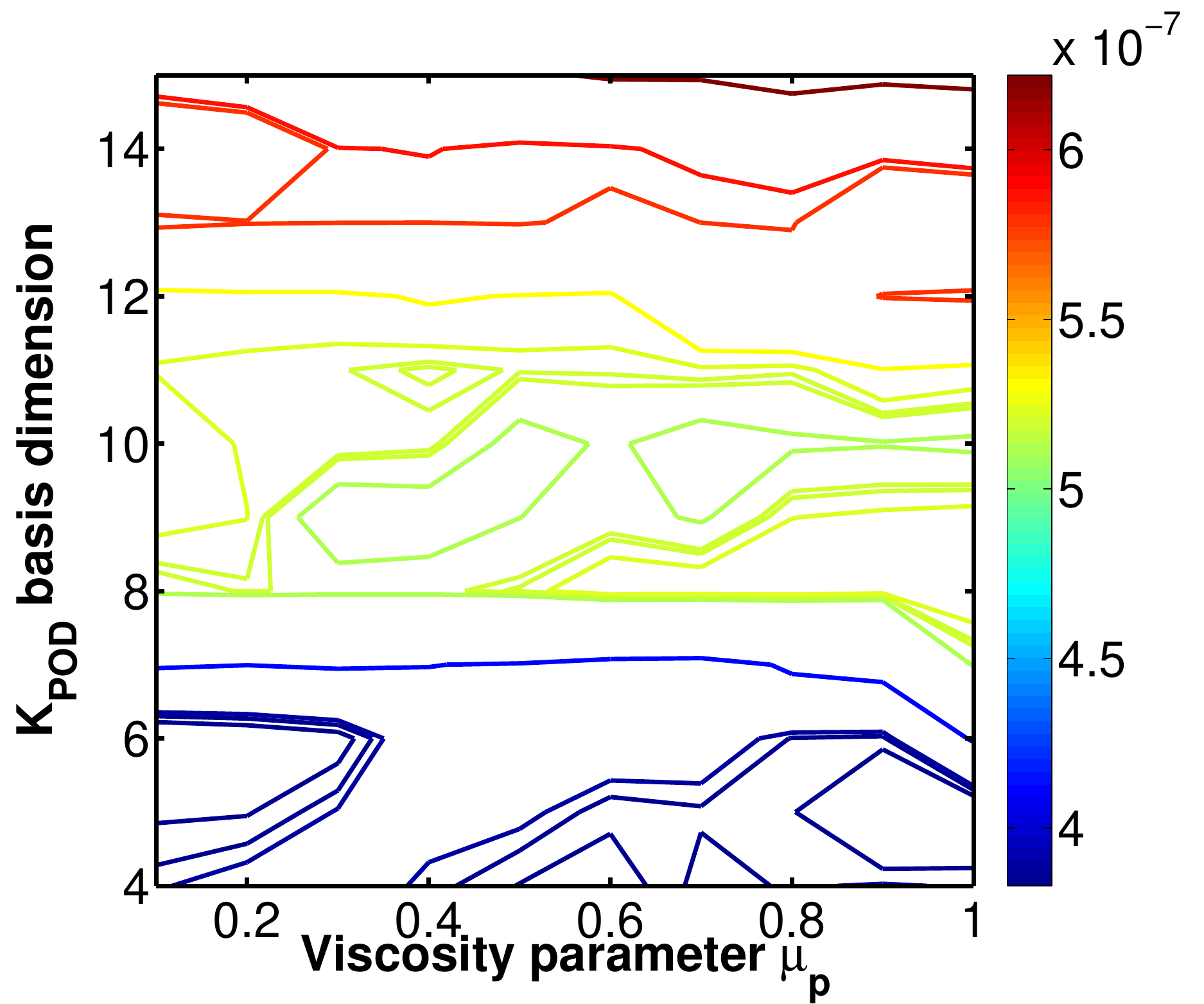}}
\caption{Isocontours for the $\textnormal{VAR}_{\rm fold}$ using ANN method.
\label{fig:VARcontours_NN}}
\end{figure}

Figures \ref{fig:contours_GP}-\ref{fig:VARcontours_NN} shows the isocontours of the $\textnormal{E}_{\rm fold}$ and $\textnormal{VAR}_{\rm fold}$ computed using \eqref{eqn:err_fold} for different $K_{POD} $ and $\mu_p$ using ROMES, MFC, and MP-LROM models constructed using GP and ANN methods. In total there are $12 \times 10$
configurations corresponding to different $K_{POD} $ and $\mu_p$ and as many ROMES and MFC models. The MP-LROM models are global in nature and the training set is the whole original data set. The testing set is the same for all the compared models and differs from the training sets. We can see that MP-LROM models are more accurate than those obtained via ROMES and MFC models. Including more samples associated with various POD basis sizes and $\mu_p$ is benefic. We also trained and tested all the models using five-fold cross-validation.  The average error and variance of all 120 $ \textnormal{E}_{\rm fold}$s and $\textnormal{VAR}_{\rm fold}$s are compared against those obtained using MP-LROM error models and are summarized  in tables  \ref{tab:experm_ROMES_mul_Efold} and \ref{tab:experm_ROMES_mul_Vfold}. This shows that the MFC models outperform the ROMES ones, for our experiment, and the MP-LROM models are the most accurate. The MP-LROM models perform better since they employ more features and samples than the other models which help the error models tune the parameters better. We also notice the efficiency of the MFC models from accuracy point of view considering that they use very few samples. In the case of large parametric domains the MP-LROM error models may require a very large data set with a lot of features. By using only subsets of the whole data set near the vicinity of the parameters of interest and applying the active subset method \cite{constantine2014active} can help prevent the potential curse of dimensionality problem that MP-LROM might suffer.

%
%
%
%
%
%
%
\begin{table}[H]
\begin{center}
    \begin{tabular}{  | c  | c | c | c |}
    \hline
   & ROMES  & MFC  & MP-LROM
      \\ \hline
 ANN & $0.3844$ & $ 0.0605$ & $ 8.8468  \times 10 ^{-4}$
     \\ \hline
 GP & $ 0.2289$ & $ 0.0865 $& $ 0.0362 $
 \\ \hline
     \end{tabular}
\end{center}
 \caption{Average error of all 120 $ \textnormal{E}_{\rm fold}$s for three methods.}
  \label{tab:experm_ROMES_mul_Efold}
\end{table}
\begin{table}[H]
\begin{center}
    \begin{tabular}{ | c  | c | c | c | }
    \hline
   & ROMES  & MFC  & MP-LROM      \\ \hline
 ANN & $0.0541$ & $ 0.0213$ & $ 4.9808  \times 10 ^{-7}$
     \\ \hline
 GP & $ 0.0051$ & $ 0.0049 $& $5.4818  \times 10 ^{-4}$
 \\ \hline
     \end{tabular}
\end{center}
 \caption{Average variance of all 120 $ \textnormal{VAR}_{\rm fold}$s for three methods}
  \label{tab:experm_ROMES_mul_Vfold}
\end{table}

\begin{figure}[h]
  \centering
  \subfigure[GP error model] {\includegraphics[scale=0.35]{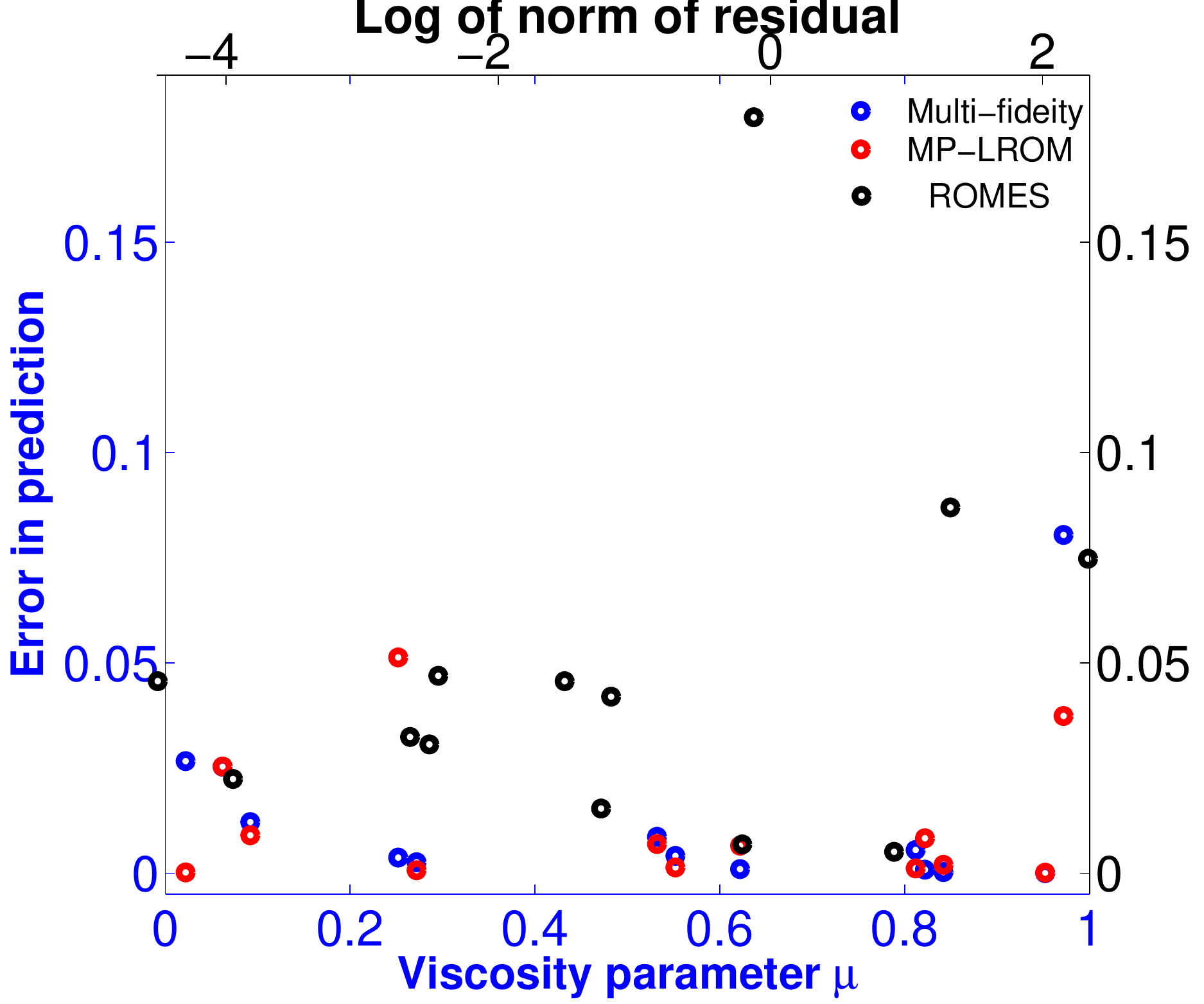}}
  \subfigure[ANN error model]{\includegraphics[scale=0.35]{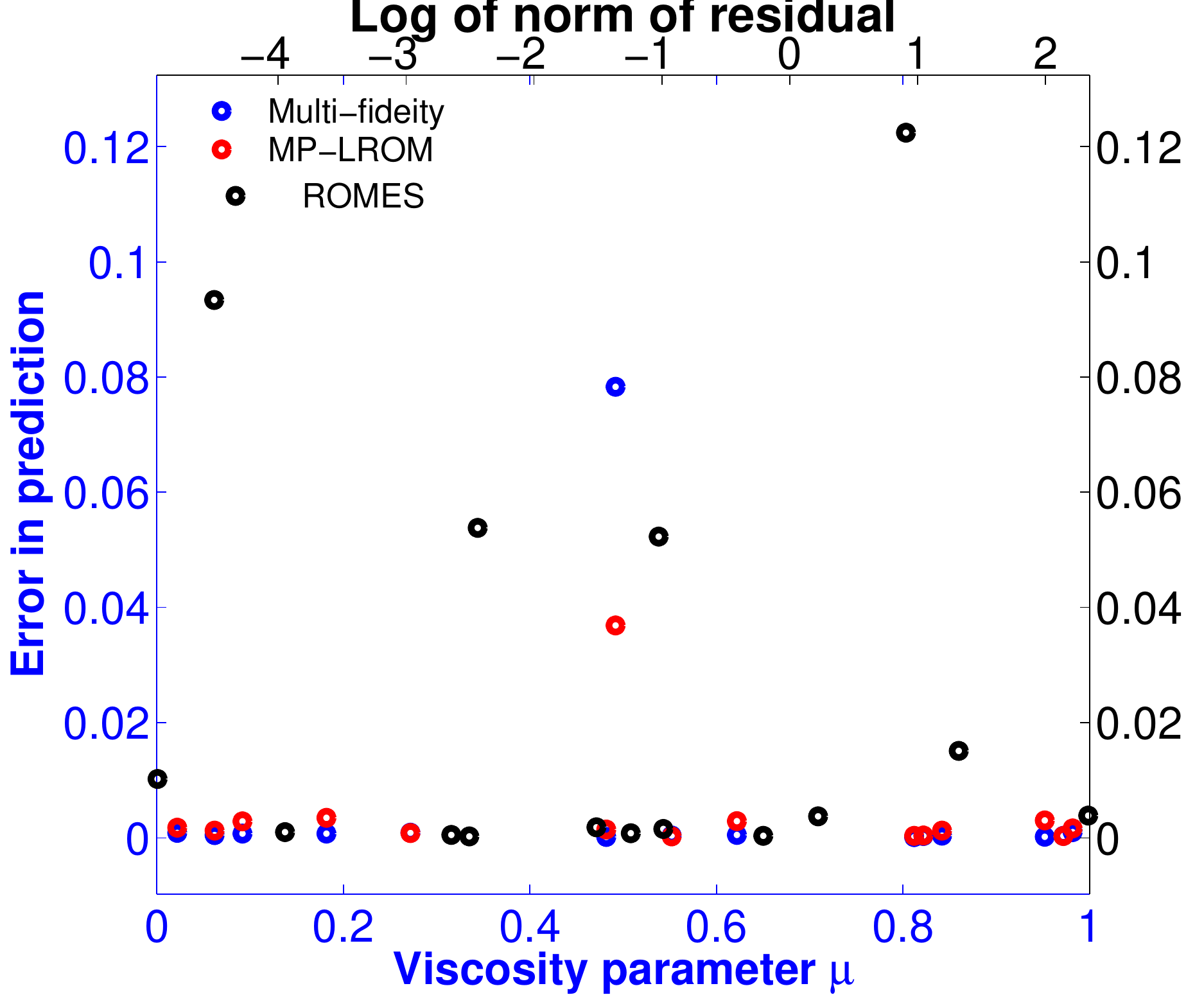}}
\caption{The error in predictions of all methods using $K_{POD} = 10$ and $\mu_p=1$.
For GP error models the overall average of errors in predictions is  $0.0131, 0.0487, 0.0095$ for MFC, ROMES and MP-LROM respectively.
For ANN error models the overall average of errors in predictions is  $0.0056, 0.0240, 0.0029$ for MFC, ROMES and MP-LROM respectively. The top x-axis shows the corresponding logarithms of the residuals norms used as inputs for the ROMES method.}
\label{fig:error_MULROMES}
\end{figure}

Finally we compare and show the average of the errors in prediction of five different errors models designed using
ROMES, MFC, and MP-LROM methods for one of the subsets corresponding to $K_{POD} = 10$ and $\mu_p=1$. The testing set is randomly selected from the samples and is not included in the training sets. The training set for both ROMES and MFC models are the same. In order to prevent the bias in prediction, each time the error models are trained on randomly selected $80\%$ of the training sets and tested with the fixed test set. We repeated this five times and the average of error in prediction is obtained. Figure \ref{fig:error_MULROMES} shows the average of error in prediction for all models implemented using GP and ANN methods.

%
Since we are including more features in our mappings, we achieve more accurate predictions compared to other existing methods such as ROMES and MFC. However, there is always a trade-off between the computational complexity and the accuracy. For more accurate results, one can generate a bigger dataset with more samples taken from the parameter domain of the underlying model. This elevates the computational complexity since the dataset requires more high-fidelity model solutions and the probabilistic mappings are more costly to construct in the training phase.
Techniques such as principal component analysis and active subspace can alleviate the curse of dimensionality for big datasets by selecting the most effective features and ignoring the less-effective ones.

%
%

\subsubsection{Selecting the dimension of reduced-order model}
\label{sect:optimal_base}
 Here we construct MP-LROM models to predict the reduced basis dimension that account for a-priori specified accuracy levels in the reduced-order model solution. The models are constructed using GP and ANN methods and have the following form
  \begin{equation}\label{eqn:MP-LROM_dimension}
  \phi_{MP-LROM}^d: \{\mu_p,\log{{\varepsilon}}_{\mu_p,\mu_p,K_{POD}}^{HF}\} \mapsto \widehat{K_{POD}}.
\end{equation}
 The input features of this model consist of the viscosity parameter $\mu_p \in [0.01,1]$ and the log of the Frobenius norm of the error between the high-fidelity and reduced-order models \eqref{eqn:param_rang_err}.
 The searched output $\widehat{K_{POD}}$ is the estimation of the dimension of the reduced manifold $K_{POD}$.
The data set contains equally distributed values of $\mu_p$ over the entire parametric domain
 $\, \mu_p \in \{0.01, 0.0113,0.0126,\ldots,0.9956 \}$, reduced basis dimensions $K_{POD}$ spanning the set $\{4,5,\ldots,14,15\}$ and the logarithm of the reduced-order model error $\log \varepsilon_{\mu_p,\mu_p,K_{POD}}^{HF}$. We use GP and ANN methods to construct two MP-LROM models to predict the dimension of local reduced-order models given a prescribed accuracy level.

\begin{figure}[h]
  \centering
  \includegraphics[width=0.5\textwidth, height=0.40\textwidth]{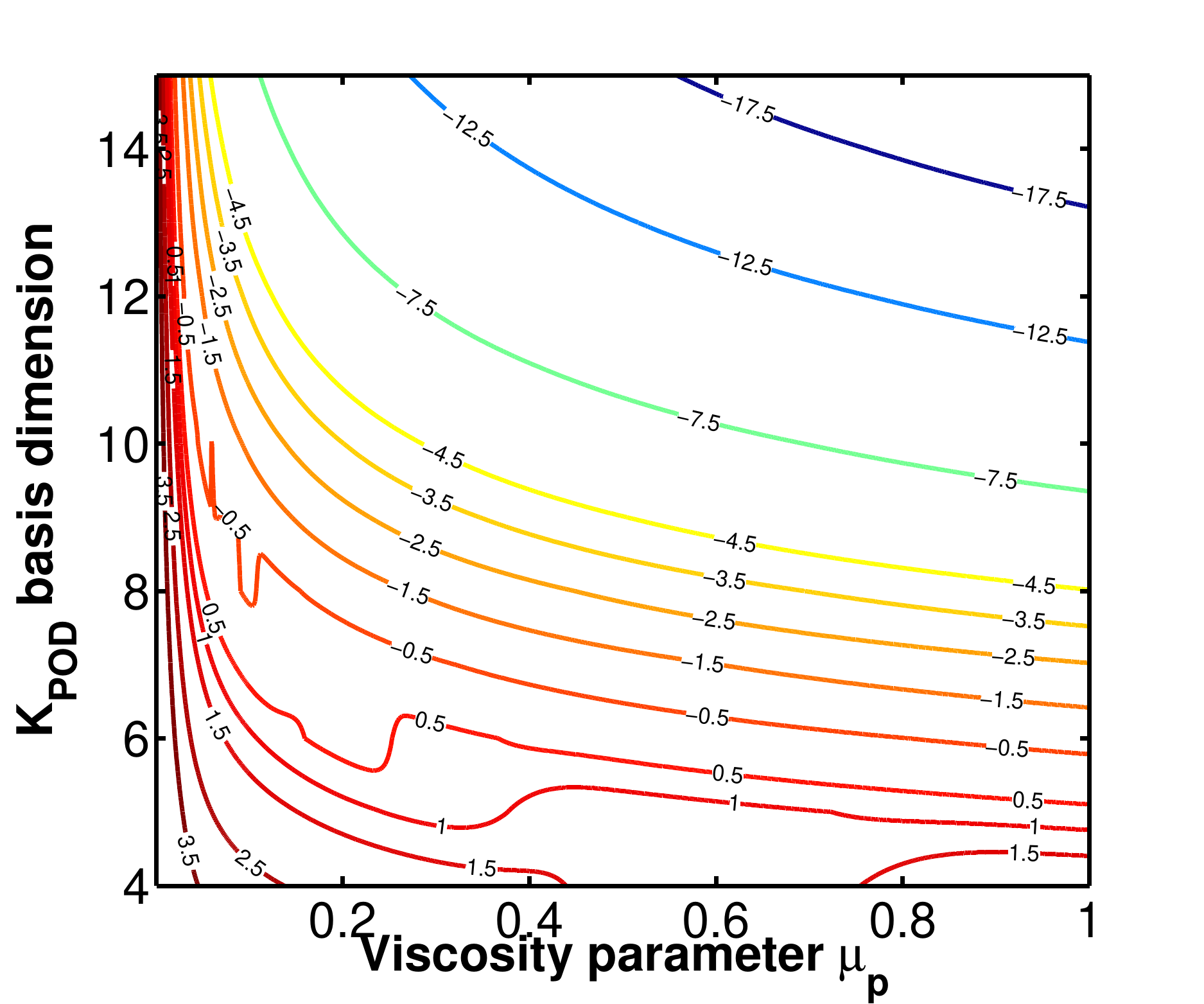}
\caption{Isocontours of the reduced model errors for different POD basis dimensions and viscosity parameters $\mu_p$. \label{fig:basis_contour_log} }
\end{figure}

During the training phase, the MP-LROM models will learn the dimensions of reduced-order basis $K_{POD}$ associated with the parameter $\mu_p$ and the corresponding error $\log\varepsilon_{\mu_p,\mu_p,K_{POD}}^{HF}$. Later they will be able to estimate the proper dimension of reduced basis by providing it the specific viscosity parameter $\mu_p$ and the desired precision $\log\bar{\varepsilon}$. The computational cost is low once the models are constructed. The output indicates the dimension of the reduced manifold for which the ROM solution satisfies the corresponding error threshold. Thus we do not need to compute the entire spectrum of the snapshots matrix in advance which for large spatial discretization meshes translates into important computational costs reduction.
 Figure \ref{fig:basis_contour_log} illustrates the contours of the log of reduced-order model errors
over all the values of the viscosity parameter $\mu_p \in \{0.01, 0.0113,0.0126\ldots 1\}$ and various POD dimensions $K_{POD} = \{4,5,\ldots,14,15\}$.

A neural network with $5$ hidden layers and hyperbolic tangent sigmoid activation function in each layer is used while for the Gaussian Process we have used the squared-exponential-covariance kernel \eqref{eq_cov}.
 For both MP-LROM models, the results were rounded such as to generate natural numbers.
Table \ref{tab:Opt_log} shows the average and variance of error in GP and ANN predictions using different sample sizes. ANN outperforms GP and as the number of data points grows, the accuracy increases and the variance decreases. The results are obtained using a conventional validation with $80\% $ of the sample size dedicated for training data and the other $20\% $ for the test data. The employed formula is described in equation \eqref{eqn:err_fold}.

%

%
\begin{table}[H]
\begin{center}
    \begin{tabular}{ | l | l | l |  l | l |}
    \hline
     & \multicolumn{2}{|c|}{MP-LROM GP} & \multicolumn{2}{|c|}{MP-LROM ANN} \\
 \hline
sample size &  $\textnormal{E}_{\rm fold}$   &  $\textnormal{VAR}_{\rm fold}$    & $\textnormal{E}_{\rm fold}$   &  $\textnormal{VAR}_{\rm fold}$
     \\ \hline
 100 & $ 0.2801 $ & $0.0901$ & $ 0.1580$ & $ 0.02204 $
     \\ \hline
 1000 & $0.1489$ & $ 0.0408 $ & $ 0.0121 $ & $ 0.0015 $
 \\ \hline
 3000 & $0.1013 $ & $ 0.0194 $ & $ 0.0273 $ & $ 0.0009 $
 \\ \hline
 5000 & $ 0.0884 $ & $ 0.0174 $ & $ 0.0080 $ & $ 0.0002 $
 \\ \hline
     \end{tabular}
\end{center}
 \caption{ Average and variance of errors in prediction of reduced basis dimension using MP-LROM  models for different sample sizes}
  \label{tab:Opt_log}
\end{table}

\begin{figure}[h]
  \centering
  \subfigure[$100$ samples] {\includegraphics[scale=0.35]
{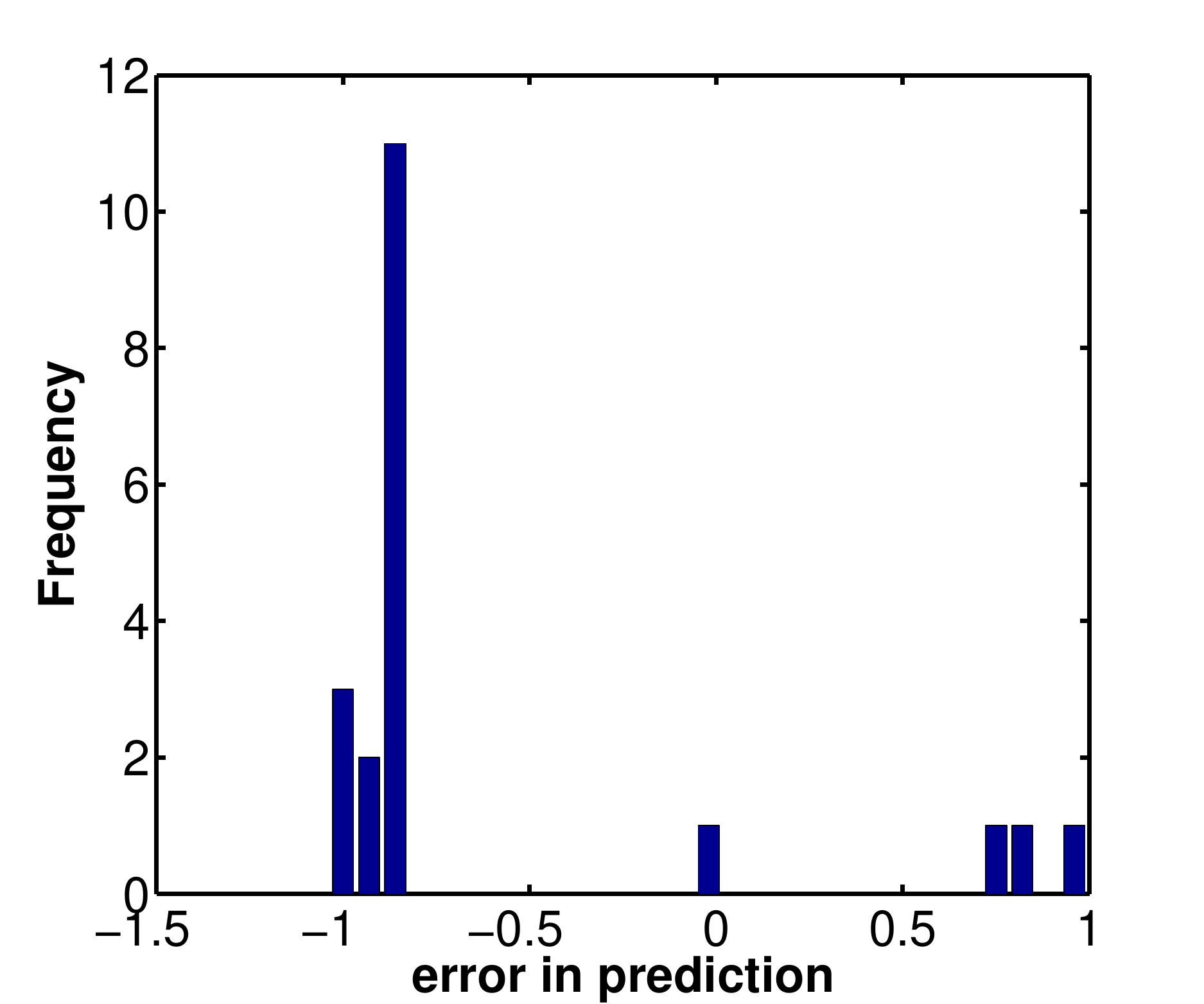}}
  \subfigure[$1000$ samples] {\includegraphics[scale=0.35]
{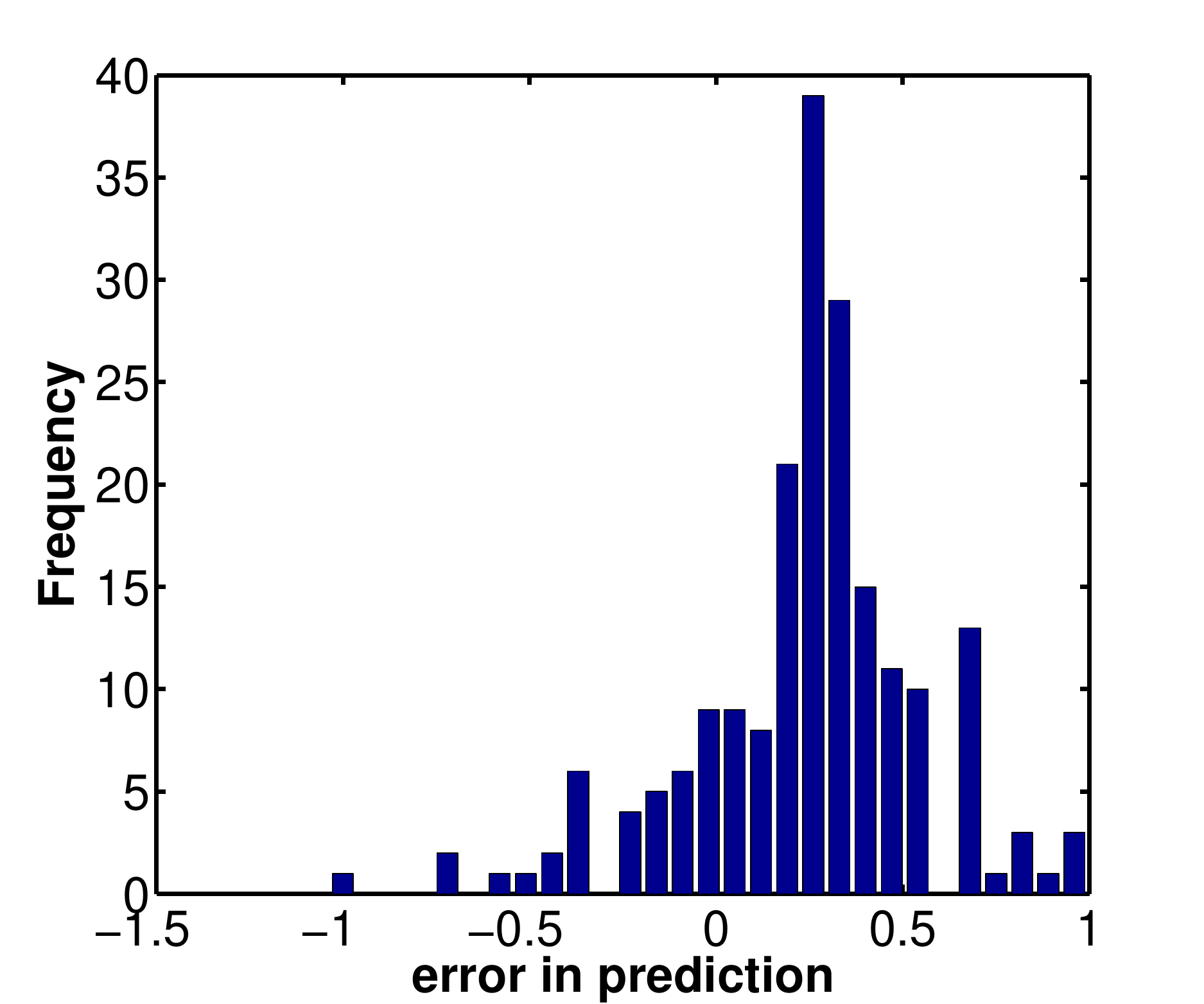}}
\caption{ Histogram of errors in prediction of the reduced basis dimension using ANN MP-LROM for different sample sizes
\label{fig:hist_NNDim}}
\end{figure}

\begin{figure}[h]
  \centering
  \subfigure[$100$ samples] {\includegraphics[scale=0.35]
{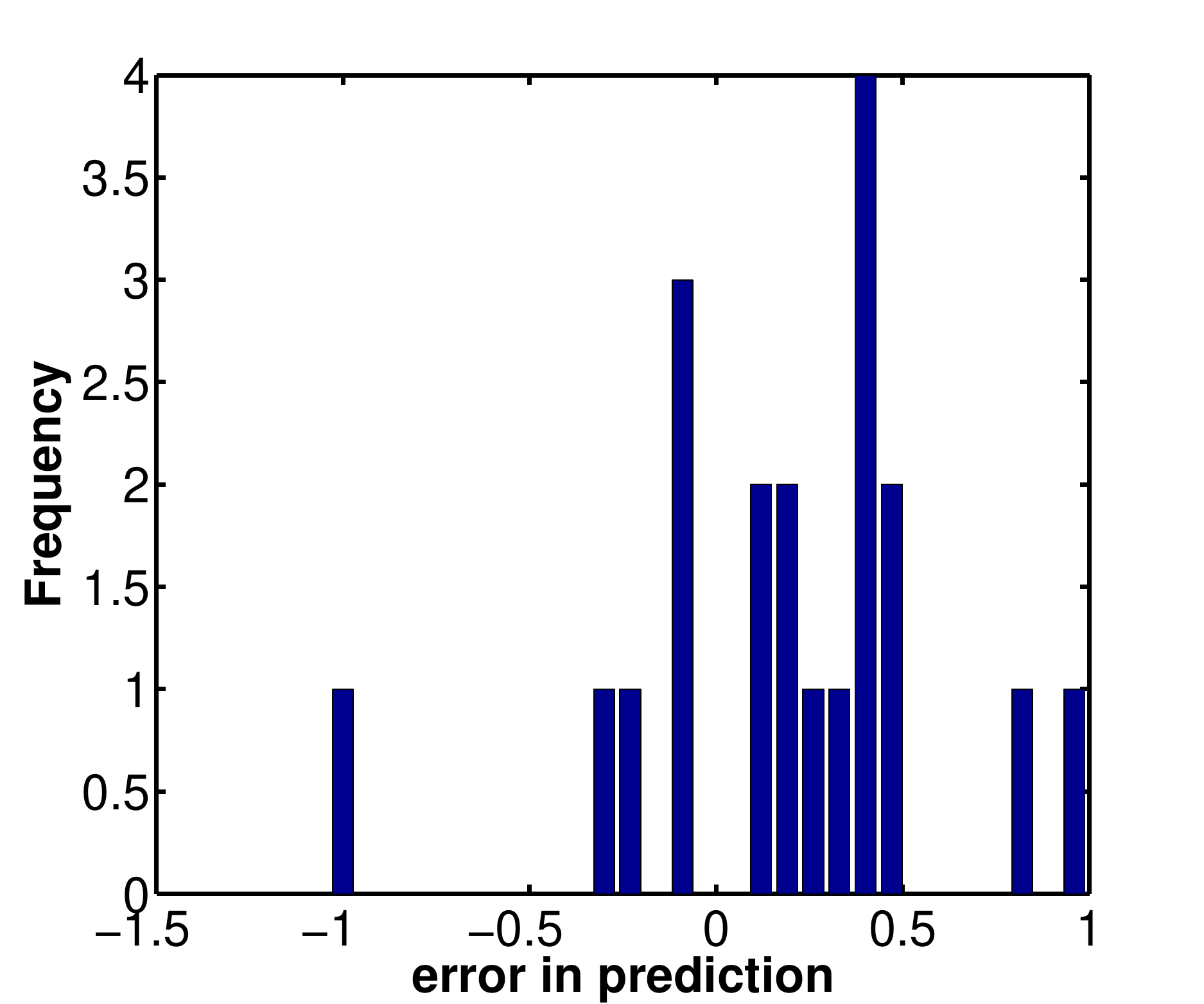}}
    \subfigure[$1000$ samples] {\includegraphics[scale=0.35]
{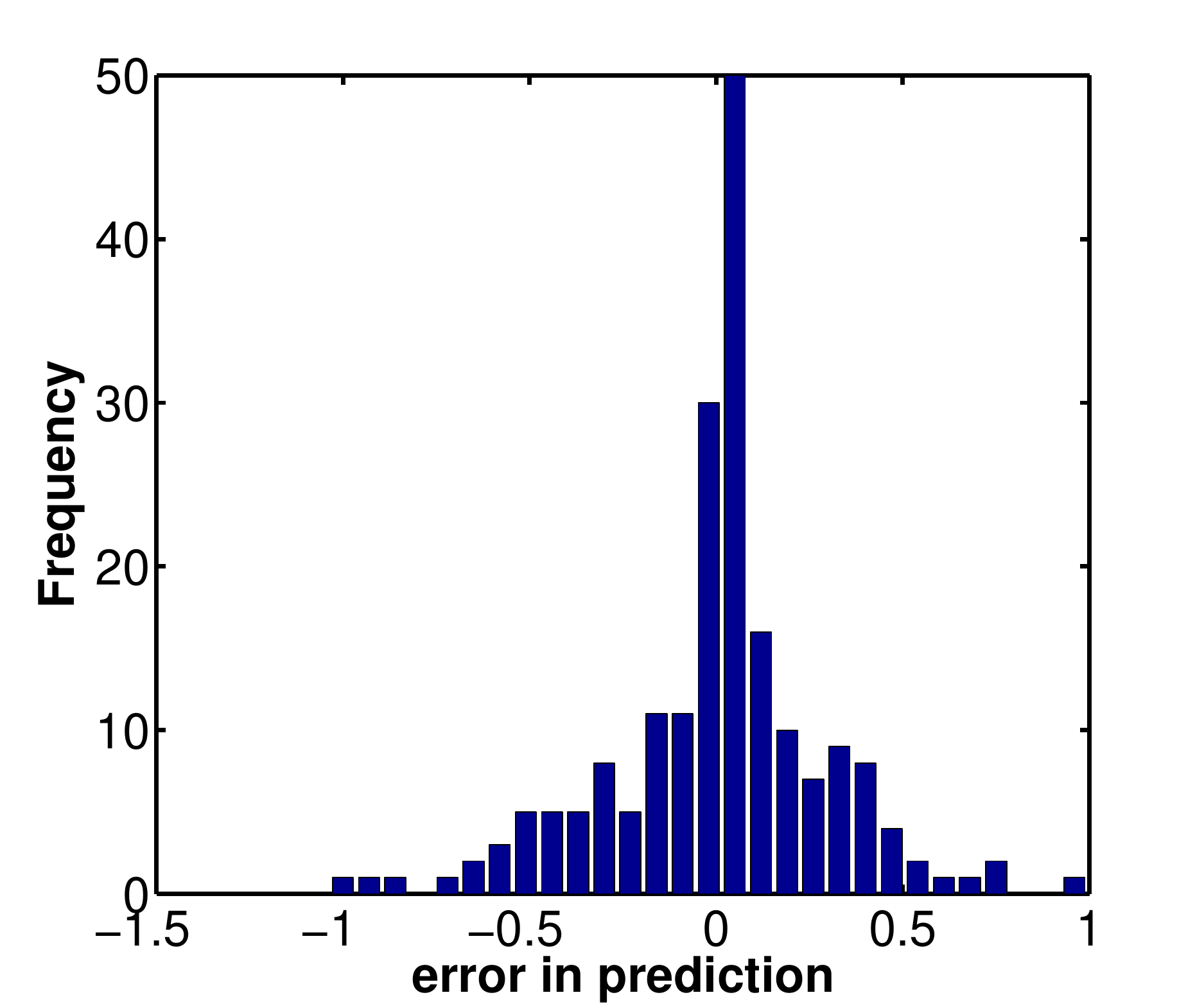}}
\caption{ Histogram of errors in prediction of the reduced basis dimension using GP ML-ROM for different sample sizes
\label{fig:hist_GPDim}}
\end{figure}

Figures \ref{fig:hist_NNDim} and \ref{fig:hist_GPDim} show the prediction errors using $100$ and $1000$ training samples for both MP-LROM models constructed via ANN and GP models. The histograms shown in Figure \ref{fig:hist_GPDim}, as stated before, can assess the validity of GP assumptions. Once the number of samples is increased, the data set distribution shape is closer to the Gaussian profile $\mathcal{N} (0, \sigma_n^2)$ than in the case of the data set distribution shown in Figure \ref{fig:ParamHist_GP} used for generation of MP-LROM models for the prediction of local reduced-order model errors.
%


%


To assess the accuracy of the MP-LROM models, the data set is randomly partitioned into five equal size sub-samples, and  five-fold cross-validation test is implemented. The five results from the folds are averaged and they are presented in Table \ref{tab:experm1}. The
 ANN model correctly estimated the dimension of the reduced manifold in $87\%$ cases.  GP correctly estimates the POD dimension $53\%$ of the times. The variance results shows that the GP model has more stable predictions indicating a higher bias in the data.

\begin{table}[H]
\begin{center}
\begin{small}
    \begin{tabular}{ | c | p{1.55cm} | p{1.25cm} | c | c | c | p{1.55cm} |  p{2.75cm} |}
    \hline
  Dimension discrepancies& zero & one  & two  &  three  & four  & $>$ four & $ VAR $
      \\ \hline
 ANN MP-LROM  & $87\% $ & $11\%$ & $2 \%$ & 0 & 0 & 0 & $2.779 \times 10^{-3}$
     \\ \hline
 GP MP-LROM  & $53\%$ & $23 \%$ & $15\%$ & $5 \%$ &$ 3\%$ & $1 \%$ & $4.575 \times 10^{-4}$
     \\ \hline
    \end{tabular}
\end{small}
\end{center}
 \caption{POD basis dimension discrepancies between the MP-LROM predictions and true values over five-fold cross-validation. The errors variance is also computed.}
  \label{tab:experm1}
\end{table}

\begin{center}
\begin{figure}[H]
	\begin{centering}
	\includegraphics[width=0.5\textwidth, height=0.4\textwidth]{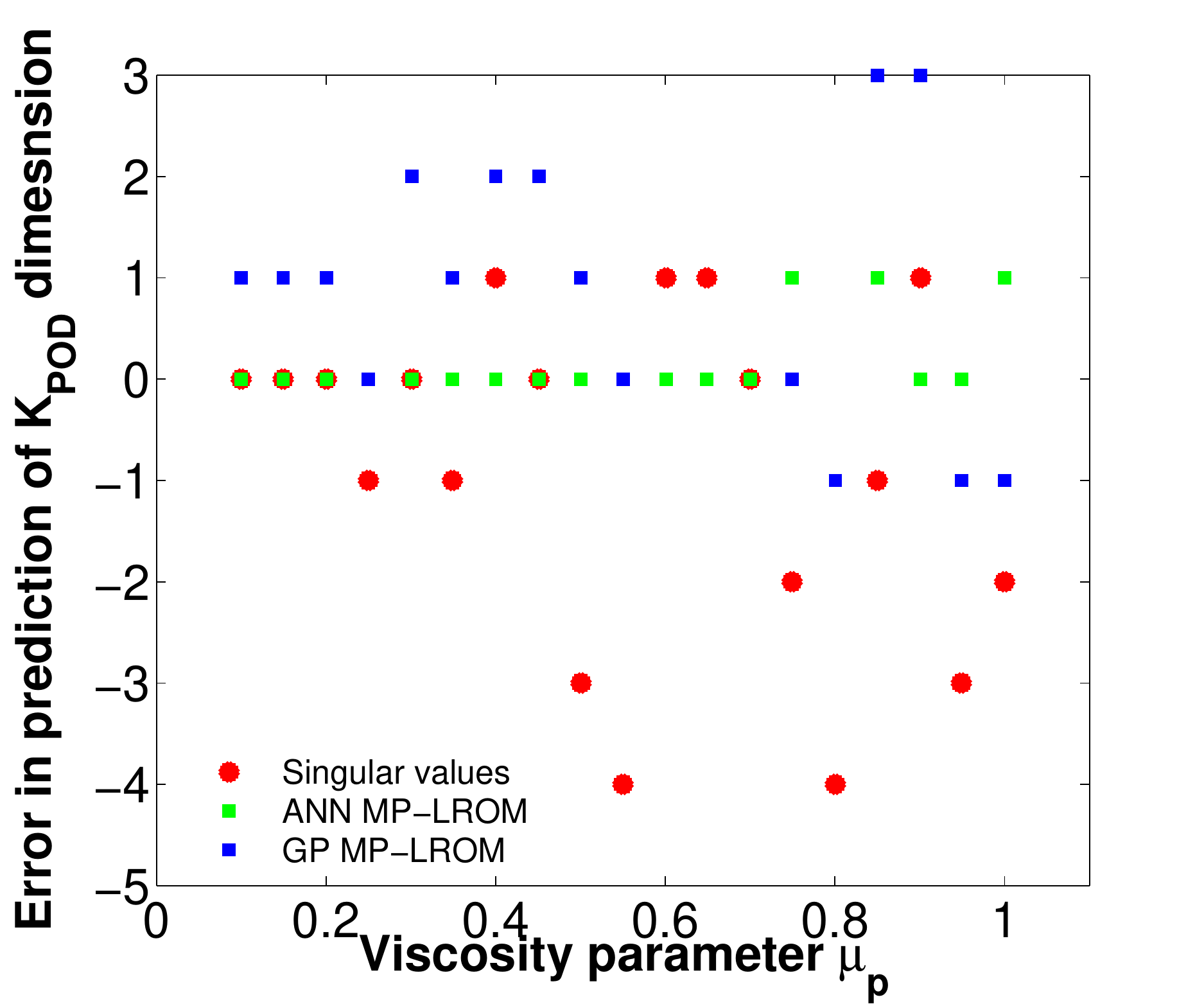}
    \caption{Average error of the POD dimension prediction on a randomly selected test data with desired accuracy of $\bar{\varepsilon}=10^{-3}$.
The average of the absolute values of the error in prediction are 1.31, 0.21 and 1.38 for singular value based method and MP-LROM models constructed using ANN and GP.}
	\label{fig:expm1_pod}
	\end{centering}
\end{figure}
\end{center}

In Figure \ref{fig:expm1_pod}, we compare the output of the MP-LROM models against the { singular values} based estimation on a set of randomly selected test data. The estimation derived from {the singular values}  is the standard method for selecting the reduced manifold dimension when a prescribed level of accuracy of the reduced solution is desired.  Here the desired accuracy $\bar{\varepsilon}$ is set to $10^{-3}$. The mismatches between the predicted and true dimensions are depicted in Figure \ref{fig:expm1_pod}.  The predicted values are the averages over five different MP-LROM models constructed using ANN and GP methods. The models were trained on random $80 \%$ split of data set and tested on the fixed selected $20 \%$ test data. We notice that the snapshots matrix spectrum underestimates the true dimension of the manifold as expected since the `in-plane' errors are not accounted. The ANN predictions were extremely accurate for most of the samples while the GP usually overestimated the reduced manifold dimensions.

\section{Conclusions}
\label{sect:conc}

In this study, we introduced new multivariate input-output models (MP-LROM) to predict the errors and dimensions of local parametric reduced-order models. Approximation of these mappings were built using Gaussian Process and Artificial Neural Networks. 
Initially, we compared our MP-LROM error models against those constructed with multi-fidelity correction technique (MFC) and reduced order model error surrogates method (ROMES). Since global bases are used by MFC and ROMES methods, we implemented corresponding local error models using only small subsets of the data utilized to generate our MP-LROM models. In contrast, the MP-LROM models are global and rely on a global database. Moreover, our MP-LROM models differ from the ROMES \cite{drohmann2015romes} and MFC models \cite{alexandrov2001approximation}, having more additional features such as reduced subspace dimension and are specially projected for accurate predictions of local parametric reduced-order models errors. As such, the MP-LROM models require significantly more and different data than MFC models. The numerical experiments revealed that our MP-LROM models are more accurate than the models constructed with MFC and ROMES methods for estimating the errors of local parametric reduced-order 1D-Burgers models with a single parameter.  In the case of large parametric domains, the MP-LROM error models could be affected by the curse of
dimensionality due to the large number of input features. In the future we plan to use only subsets of the global data set near the vicinity of the parameters of interest and combine our technique with the active subspace method \cite{constantine2014active} to prevent the potential curse of dimensionality that the MP-LROM models might suffer.

Next we addressed the problem of selecting the dimension of a local reduced-order model when its solution must satisfy a desired level of accuracy. The approximated MP-LROM models based on Artificial Neural Networks better estimated the ROM basis dimension in comparison with the results obtained by truncating the spectrum of the snapshots matrix.

In the future we seek to decrease the computational complexity of the MP-LROM error models. Currently the training data required by the machine learning regression MP-LROM models rely on many high-fidelity simulations. By employing error bounds, residual norms \cite{drohmann2015romes} and a-posteriori error estimation results \cite{Volwein_aposteriori_2016,nguyen2009reduced}, this dependency could be much decreased. On-going work focuses on applications of MP-LROM error model. We are currently developing several algorithms and techniques that employ MP-LROM error model as a key component to generate decomposition maps of the parametric space associated with accurate local reduced-order models.


 In addition, we plan to construct machine learning MP-LROM models to estimate the errors in quantities of interest computed with reduced-order models. The predictions of such error models can then be used to speed up the current trust-region reduced-order framework \cite{Arian_2000,bergmann2008optimal} by eliminating the need of high-fidelity simulations for the quality evaluation of the updated controls.

\section*{Acknowledgements}
This work was supported in part and by the award NSF CCF 1218454 and by the Computational Science Laboratory at Virginia Tech.

\label{sect:bib}
\bibliographystyle{plain}

\bibliography{ML_bib,Additions_sandu,CDS_E_proposal,comprehensive_bibliography1,data_assim_fdvar,data_assim_weak-fdvar,NSF_KB,POD_bib,Razvan_bib,Razvan2_bib,Razvan_bib_ROM_IP,Razvan_update_bib,reduced_models,ROM_state_of_the_art,sandu,Software}

\end{document}